\documentclass[11pt,twoside]{amsart}
\usepackage{amssymb,amsmath,amscd,amsthm,latexsym,mathtools,mathrsfs}
\usepackage[noadjust]{cite}
\usepackage{indentfirst,graphicx,color,enumerate,array}
\numberwithin{equation}{section}
\theoremstyle{plain}
\newtheorem{theorem}{Theorem}[section]
\newtheorem{lemma}[theorem]{Lemma}
\newtheorem{proposition}[theorem]{Proposition}

\theoremstyle{definition}
\newtheorem{definition}[theorem]{Definition}
\newtheorem{remark}[theorem]{Remark}
\usepackage[hidelinks]{hyperref}

\newcommand{\R}{\mathbb R}
\newcommand{\C}{\mathbb C}
\newcommand{\norm}[1]{\lVert #1\rVert}
\newcommand{\abs}[1]{\left|#1\right|}
\newcommand{\dd}{\,\mathrm d}
\newcommand{\D}{\mathcal D}

\title[Impact Mechanism of a Point Interaction in the Cubic Half-Wave Equation]
{Impact Mechanism of a Point Interaction in the Cubic Half-Wave Equation: Well-Posedness, Constrained Energy Levels and Near-Critical Rigidity}

\author[X. Ju]{Xingkai Ju}
\address{School of Mathematics, Jilin University, Changchun 130012, China}
\email{xingkaiju\_math@163.com}

\author[Y. Wei]{Yuanhong Wei}
\address{School of Mathematics \& International Institute for Mathematical Sciences,\hfill\break
Jilin University, Changchun 130012, China}
\email{weiyuanhong@jlu.edu.cn}

\begin{document}

\begin{abstract}
We study the one-dimensional nonlinear half-wave equation with a
logarithmically renormalized point interaction. The nonlinearity is focusing and cubic, which is mass-critical. We establish well-posedness in
the energy space using a new method for unconditional uniqueness that requires
no auxiliary Strichartz assumption. We reveal the mechanism by which the point
interaction affects ground states. The free critical mass is the sharp
threshold: fixed-mass ground states exist and are orbitally stable below it,
whereas at and above it the constrained energy is unbounded from below. We also
prove a near-critical rigidity result showing that bounded energy concentrating
sequences converge to free ground states after rescaling.
\end{abstract}

\maketitle



\tableofcontents
\newpage

\section{Introduction}

We study the Cauchy problem for the one-dimensional cubic half-wave equation with a point interaction
\begin{equation}\label{eq:evolution}
 \left\{
 \begin{array}{@{}l@{}}
  i\partial_t u=|D|_\alpha u-|u|^{2}u,\\
  u(0,x)=u_0(x),
 \end{array}
 \right.
 \qquad (t,x)\in\R\times\R.
\end{equation}
Here \(u=u(t,x):\R\times\R\to\C\) is the complex-valued wave amplitude
and \(u_0\) is the
initial datum.  The term \(-\abs u^2u\) is the focusing cubic Kerr nonlinearity.
The operator \(i\partial_t-|D|_\alpha\) is the half-wave operator with a point
interaction at the origin, and \(\alpha\in\R\) represents the interaction strength.

In the absence of the point interaction, the linear part is the free operator \(|D|\).  This is the classical half-Laplacian, namely
\(|D|=(-\Delta)^{1/2}\). From the Fourier transform convention
\begin{equation}\label{eq:fourier-convention}
 \widehat f(\xi):=\int_\R e^{-ix\xi}f(x)\,\dd x,
 \qquad
 f(x)=\frac1{2\pi}\int_\R e^{ix\xi}\widehat f(\xi)\,\dd\xi,
\end{equation}
it holds that
\begin{equation}\label{eq:intro-free-half-wave-operator}
 \widehat{|D|f}(\xi)=|\xi|\widehat f(\xi),
 \qquad \mathrm{Dom}(|D|)=H^1(\R),
\end{equation}
and the domain of its quadratic form is \(H^{1/2}(\R)\).

One may expect that a point interaction perturbation of a self-adjoint operator
\(\mathcal A_0\) can be written formally as \(\mathcal A_0+\alpha\delta_0\), where
\(\delta_0\) is the delta function at the origin.  This description applies to
some standard cases, such as the one-dimensional Laplacian with a point
interaction.  However, it is misleading here, since a form term involving
\(u(0)\) is not well-defined on \(H^{1/2}(\R)\).  In fact, elements of
\(H^s(\R)\) with \(s>1/2\) are continuous by the Sobolev embedding.  At the
critical regularity \(s=1/2\), point evaluation is no longer continuous.
Instead, we consider the closure of the symmetric restriction of \(|D|\) to
\(C_c^\infty(\R\setminus\{0\})\).  Its self-adjoint extensions form a family
of point interaction operators.  The extension parameter is fixed by a
renormalized boundary condition at the origin.

We explain the definition of the linear operator
\(|D|_\alpha\) in detail. Let us consider the free Green function
\begin{equation}\label{eq:intro-half-wave-green}
 G_\lambda=(|D|+\lambda)^{-1}\delta_0,
 \qquad
 G_\lambda(x)=\frac1{2\pi}\int_{-\infty}^{+\infty}
 \frac{e^{ix\xi}}{|\xi|+\lambda}\,\dd\xi,
 \qquad x\neq0,\quad \lambda>0.
\end{equation}
The Fourier integral in \eqref{eq:intro-half-wave-green} is understood as an
improper oscillatory integral.
The Green function belongs to \(L^2(\R)\), but not to \(H^{1/2}(\R)\). Besides,  we have
\begin{equation}\label{eq:intro-green-log}
 G_\lambda(x)=\frac1\pi\log\frac1{\lambda|x|}+O_\lambda(1)
 \qquad\text{as }x\to0.
\end{equation}
The logarithmic singularity makes the diagonal value of the free resolvent
divergent.  After the divergent part is absorbed into the extension parameter,
the finite remainder is encoded by the Krein principal function
\begin{equation}\label{eq:intro-principal-function}
 \Phi_\alpha(-\lambda)
 =\alpha+\frac1\pi\log \lambda
 =\frac1\pi\log\frac{\lambda}{\lambda_\alpha},
\end{equation}
where \(\lambda_\alpha:=e^{-\pi\alpha}\).  The operator \(|D|_\alpha\) has the
simple negative eigenvalue \(-\lambda_\alpha\) with eigenfunction
\(G_{\lambda_\alpha}\).
For any fixed \(\lambda>0\), write
\(u=\phi+qG_\lambda\), where \(\phi\) is the regular part and \(qG_\lambda\)
is the singular part.  The operator domain is
\begin{equation}\label{eq:intro-renormalized-operator-domain}
 \mathrm{Dom}(|D|_\alpha)
 =\{u=\phi+qG_\lambda:\ \phi\in H^1(\R),\ q\in\C,\
   \phi(0)=\Phi_\alpha(-\lambda)q\},
\end{equation}
and the action of \(|D|_\alpha\) on this domain is defined by
\begin{equation}\label{eq:intro-renormalized-operator-action}
 |D|_\alpha u=|D|\phi-\lambda qG_\lambda.
\end{equation}
The condition \(\phi(0)=\Phi_\alpha(-\lambda)q\) in
\eqref{eq:intro-renormalized-operator-domain} is the renormalized boundary
condition at the interaction point.  It involves the regular value
\(\phi(0)\), rather than \(u(0)\), since the full function \(u\) has no finite
value at the origin when \(q\neq0\).

The half-wave equation is a prototypical semirelativistic fractional nonlinear
evolution equation.  Physically, it models massless semirelativistic particles
and describes the propagation of self-focusing Kerr solitons in fractional
optics.  Mathematically, the half-wave equation lacks the standard dispersive
estimates available for the Schr\"odinger equation.  For the free
one-dimensional half-wave problem, Frank and Lenzmann proved uniqueness of the positive ground-state
profile \cite{FrankLenzmann}; Krieger, Lenzmann and Rapha\"el obtained the
cubic critical mass and minimal-mass nondispersive dynamics
\cite{KriegerLenzmannRaphael}; and the low-regularity obstruction was
established in \cite{ChoffrutPocovnicu}.  More generally, the local and global
Cauchy theory and the stability or instability of ground states were studied
in \cite{BellazziniGeorgievVisciglia}.

In Laskin's fractional Schr\"odinger framework \cite{LaskinFSE}, \(|D|\)
corresponds to L\'evy index one, and the cubic term represents the self-focusing
Kerr response \cite{ChenEtAlCriticalCollapse}.  The cubic half-wave equation
with a local defect has been studied directly in fractional optics, where an
attractive defect can stabilize pinned solitons
\cite{MayteevarunyooMalomed}.  The same
linear pair of \(|D|\) and a delta interaction also arises for massless
semirelativistic particles \cite{ErmanGadellaUncu}.  The corresponding contact interaction requires
logarithmic renormalization through dimensional transmutation
\cite{AlHashimiShalabyWiese}.  For further results on point interactions, we
refer to the classical monograph of Albeverio, Gesztesy, H{\o}egh-Krohn and
Holden \cite{AGHH}, renormalized relativistic delta Hamiltonians
\cite{AlHashimiShalabyWiese,ErmanGadellaUncu}, the classification of fractional
singular perturbations \cite{MichelangeliOttoliniScandone} and their
approximation by shrinking potentials \cite{MichelangeliScandone}.

For one-dimensional NLS equations with point interactions, Fukuizumi, Ohta
and Ozawa studied well-posedness and stability \cite{FukuizumiOhtaOzawa}.
Further results concern dynamics \cite{AdamiNojaDynamics}, ground states and
stability \cite{AdamiNojaVisciglia,AdamiBoniNakamuraRuighi} and mass-critical
blow-up \cite{TangXu}.  These models use ordinary
\(H^1\)-type boundary traces and require no logarithmic renormalization.  A
closer analogue of the present model is the two-dimensional Laplacian with a
point interaction.  Cacciapuoti, Finco and Noja established well-posedness
for the nonlinear problem \cite{CacciapuotiFincoNoja}.  Further results concern
ground states and standing waves
 \cite{AdamiBoniCarloneTentarelli,FukayaGeorgievIkeda,
 GeorgievMichelangeliScandone,BoniGallone}, and blow-up dynamics
 \cite{FincoNoja,ForcellaGeorgiev,DuanZhang}.  Related fractional models treat nonlinear
singular perturbations \cite{CarloneFincoTentarelli} and singular Hardy
potentials \cite{MizutaniYao}.

The present study is motivated by cubic half-wave equations with local defects
proposed in fractional optics.  The operator theory of the renormalized
half-wave point interaction is well understood.  Known nonlinear work is
limited mainly to numerical analysis and variational approximation
\cite{MayteevarunyooMalomed}.  However, the theoretical analysis of
nonlinear half-wave equations with a point interaction does not appear.
Accordingly, we consider the nonlinear problem in the renormalized point
interaction framework.

We now introduce the notation needed to state the main results.  The natural
energy space is the form domain, which for any
\(\lambda>0\) is given by
\begin{equation}\label{eq:intro-renormalized-form-domain}
 X_\alpha
 =H^{1/2}(\R)\dotplus\operatorname{span}\{G_\lambda\}.
\end{equation}
We denote its anti-dual by \(X_\alpha^*\).
Every \(u\in X_\alpha\) has a unique regular--singular decomposition
\(u=\phi+qG_\lambda\).  The number \(q\) is the charge associated with the point interaction.
With respect to this decomposition, the quadratic form is given by
\begin{equation}\label{eq:critical-form}
 \mathfrak t_\alpha[u]
 =\norm{|D|^{1/2}\phi}_2^2+\lambda\norm\phi_2^2-\lambda\norm u_2^2
   +\Phi_\alpha(-\lambda)|q|^2,
 \qquad u\in X_\alpha.
\end{equation}
The decomposition depends on \(\lambda\), whereas \(q\) and
\(\mathfrak t_\alpha[u]\) are independent of this choice.  For any
\(c>\lambda_\alpha\), the shifted form
\(\mathfrak t_\alpha[u]+c\norm u_2^2\) defines an equivalent norm on
\(X_\alpha\).  We fix \(c=2\lambda_\alpha\) and set
\begin{equation}\label{eq:intro-form-norm}
 \norm u_{X_\alpha}^2
 :=\mathfrak t_\alpha[u]+2\lambda_\alpha\norm u_2^2.
\end{equation}
The mass and energy of \(u\in X_\alpha\) are defined by
\begin{equation}\label{eq:intro-mass}
 M(u):=\norm u_2^2.
\end{equation}
\begin{equation}\label{eq:cubic-half-wave-energy}
 E_\alpha(u):=\frac12\mathfrak t_\alpha[u]-\frac14\norm u_4^4.
\end{equation}
For \(\rho>0\), let
\[
(\D_\rho f)(x):=\rho^{1/2}f(\rho x)
\]
denote the mass-preserving dilation.

Since the free operator \(|D|\) can be regarded as the limit of \(|D|_\alpha\) as \(\alpha\to\infty\), we set
\(|D|_\infty:=|D|\).  For \(v\in H^{1/2}(\R)\), define the free energy and its
fixed-mass level by
\[
E_\infty(v)
:=\frac12\norm{|D|^{1/2}v}_2^2-\frac14\norm v_4^4,
\qquad
I_\infty(\mu)
:=\inf_{\substack{v\in H^{1/2}(\R)\\ M(v)=\mu}}E_\infty(v).
\]
Let us consider the free half-wave equation
\begin{equation}\label{eq:intro-free-half-wave-equation}
 i\partial_t u=|D| u-|u|^{2}u,\qquad (t,x)\in\R\times\R.
\end{equation}
Then by setting \(u(t,x)=e^{i\omega t}Q(x)\) (\(\omega>0\)), we obtain
\begin{equation}\label{eq:intro-free-standing-wave}
 (|D|+\omega)Q=Q^3.
\end{equation}
By the above mass-preserving dilation we obtain
\begin{equation}\label{eq:intro-free-ground-state}
 (|D|+1)Q_*=Q_*^3.
\end{equation}
For this equation, Frank and Lenzmann proved the existence and uniqueness of a positive and even ground state \cite{FrankLenzmann}.  We still denote the unique positive and even ground state of \eqref{eq:intro-free-ground-state} by \(Q_*\). Define
\begin{equation}\label{eq:intro-free-critical-mass}
 \mu_*:=M(Q_*).
\end{equation}
Consequently, for every \(\omega>0\), \(\mu_*\) is the mass of ground states for \eqref{eq:intro-free-standing-wave}.

We first establish the well-posedness of \eqref{eq:evolution} in the natural
energy space \(X_\alpha\).

\begin{theorem}[Well-posedness]
\label{thm:intro-cubic-flow}
For every \(u_0\in X_\alpha\), there exist
\(T_\pm\in(0,\infty]\) such that the Cauchy problem
\eqref{eq:evolution} has a unique maximal solution
\[
 u\in C((-T_-,T_+);X_\alpha)
 \cap C^1((-T_-,T_+);X_\alpha^*).
\]
Moreover, the following properties hold.
\begin{enumerate}[\upshape(i)]
\item If \(u_{0,n}\to u_0\) in \(X_\alpha\) and \(u_n\) are the
corresponding maximal solutions, then every compact
\(J\subset(-T_-,T_+)\) lies in their lifespans for large \(n\), and
\(u_n\to u\) in \(C(J;X_\alpha)\).
\item Mass and energy are conserved:
\[
 M(u(t))=M(u_0),\qquad E_\alpha(u(t))=E_\alpha(u_0).
\]
\item If \(T_\pm<\infty\), then
\[
 \lim_{t\to\pm T_\pm}\norm{u(t)}_{X_\alpha}=\infty.
\]
\item If \(M(u_0)<\mu_*\), the solution is global and bounded in
\(X_\alpha\):
\[
 T_\pm=\infty
 \quad\text{and}\quad
 \sup_{t\in\R}\norm{u(t)}_{X_\alpha}<\infty.
\]
\end{enumerate}
\end{theorem}

\begin{remark}
The key point of the well-posedness result is unconditional uniqueness in the
energy class, with no auxiliary Strichartz or smoothing assumption.  The
structure of the argument suggests that it can be adapted to mass-subcritical
powers and to related models with logarithmically renormalized point
interactions.
\end{remark}

\begin{definition}[Orbital stability]\label{def:orbital-stability}
For a nonempty subset \(\mathcal S\) of \(X_\alpha\), set
\begin{equation}\label{eq:orbital-distance}
 \operatorname{dist}_{X_\alpha}(u,\mathcal S)
 :=\inf_{v\in\mathcal S}\norm{u-v}_{X_\alpha},
 \qquad u\in X_\alpha.
\end{equation}
We say that \(\mathcal S\) is orbitally stable in \(X_\alpha\) if, for every
\(\varepsilon>0\), there exists \(\delta>0\) such that, whenever
\(u_0\in X_\alpha\) satisfies
\[
 \operatorname{dist}_{X_\alpha}(u_0,\mathcal S)<\delta,
\]
the corresponding solution \(u\) of \eqref{eq:evolution} is global and satisfies
\begin{equation}\label{eq:orbital-stability}
 \sup_{t\in\R}\operatorname{dist}_{X_\alpha}
 \bigl(u(t),\mathcal S\bigr)<\varepsilon.
\end{equation}
\end{definition}

For \(\mu>0\), set
\begin{equation}
 I_\alpha(\mu)
 :=\inf\{E_\alpha(u):u\in X_\alpha,\ M(u)=\mu\}.
 \label{eq:cubic-half-wave-level}
\end{equation}
\begin{equation}
 \mathcal M_\alpha(\mu)
 :=\{Q\in X_\alpha:M(Q)=\mu,\ E_\alpha(Q)=I_\alpha(\mu)\}.
 \label{eq:cubic-half-wave-ground-set}
\end{equation}

We then characterize the critical fixed-mass variational problem and identify
the sharp threshold at the free critical mass \(\mu_*\).

\begin{theorem}[Fixed-mass ground states]
\label{thm:intro-cubic-ground-states}
For every \(0<\mu<\mu_*\), the set \(\mathcal M_\alpha(\mu)\) is nonempty
and orbitally stable in \(X_\alpha\), whereas for \(\mu\geq\mu_*\),
\(I_\alpha(\mu)=-\infty\) and \(\mathcal M_\alpha(\mu)=\varnothing\).
\end{theorem}

The following theorem gives a rigidity refinement at the critical threshold.

\begin{theorem}[Near-critical rigidity]
\label{thm:critical-rigidity}
Let \(u_n\in X_\alpha\) and suppose
\begin{equation}\label{eq:main-assumptions}
 M(u_n)\geq\mu_*,
 \qquad
 M(u_n)\to\mu_*,
 \qquad
 \sup_n\abs{E_\alpha(u_n)}<\infty,
 \qquad
 \norm{u_n}_{X_\alpha}\to\infty.
\end{equation}
Then there exist concentration scales \(\rho_n>0\), with \(\rho_n\to0\),
such that
\begin{equation}\label{eq:full-X}
 \inf_{\theta,y\in\R}
 \norm{e^{-i\theta}\D_{\rho_n}u_n-Q_*(\,\cdot-y)}_{X_\alpha}
 \to0.
\end{equation}
\end{theorem}

Based on our main results, an interesting observation shows the differences and
connections of solutions' behavior between the half-wave equation with point
interaction and the classical free half-wave equation, which further
demonstrates the intrinsic mechanism of point interactions in nonlinear quantum
mechanical models.  For well-posedness, although the
uniqueness mechanism differs from that in the free case, the critical
well-posedness framework is preserved.  For fixed-mass ground states, the
critical mass remains unchanged, while the constrained energy levels and the
ground state structure differ significantly from those in the free case.
Specifically, below the critical mass, the constrained infimum is not attained
in the free case but is attained in the point interaction problem.  At the
critical mass, it is attained in the free case, whereas it equals
\(-\infty\) in the point interaction problem; see
Table~\ref{tab:constrained-energy-levels}.  Near the critical
mass, concentrating sequences recover free ground states after rescaling.
\begin{table}[ht]
\centering
\small
\renewcommand{\arraystretch}{1.5}
\begin{tabular}{|c|c|c|}
\hline
\textbf{Mass range}
& \textbf{Free level \(I_\infty(\mu)\)}
& \textbf{Point interaction level \(I_\alpha(\mu)\)} \\
\hline
\(0<\mu<\mu_*\)
& \(I_\infty(\mu)=0\), not attained
& \(-\infty<I_\alpha(\mu)<0\), attained \\
\(\mu=\mu_*\)
& \(I_\infty(\mu)=0\), attained
& \(I_\alpha(\mu)=-\infty\), not attained \\
\(\mu>\mu_*\)
& \(I_\infty(\mu)=-\infty\), not attained
& \(I_\alpha(\mu)=-\infty\), not attained \\
\hline
\end{tabular}\vspace{0.2cm}
\caption{Constrained energy levels for the free half-wave and point interaction problems.}
\label{tab:constrained-energy-levels}
\end{table}

Taken together, the main results of this paper establish a systematic theory for the
cubic half-wave equation with a point interaction, involving
well-posedness, fixed-mass ground state theory and near-critical rigidity.
The principal new ingredients are unconditional uniqueness in the natural energy
space and a new Green charge analysis for the critical and near-critical
regimes.  Unconditional uniqueness asserts that any two solutions in the natural energy space emanating from the same initial datum coincide, without any auxiliary Strichartz or other smoothing assumptions.
The main difficulties stem from the logarithmically renormalized point
interaction, the absence of standard Strichartz estimates and the mass-critical
cubic nonlinearity.

We observe that existing well-posedness results for nonlinear Schr\"odinger equations with a renormalized linear point interaction rely on Strichartz estimates. For the half-wave equation with a point interaction, such Strichartz estimates are inapplicable, which constitutes the principal challenge for establishing well-posedness. We establish unconditional uniqueness directly in the natural energy space without invoking Strichartz estimates at any stage of the proof.  For the free cubic
half-wave equation, Krieger, Lenzmann and Rapha\"el also
proved unconditional uniqueness in \(H^{1/2}(\R)\) without using Strichartz
estimates \cite{KriegerLenzmannRaphael}.  Their proof employs a Yudovich argument based on finite exponent Sobolev bounds. However, it cannot be extended to our setting, since the logarithmic Green component produces a squared logarithmic loss and the corresponding Yudovich estimate cannot be closed.  We instead remove the bound
state pole and prove a localized resolvent estimate on the continuous spectrum,
which gives Kato smoothing for the singular component.  The Green function
decomposition reduces the remaining loss to a single logarithm, and the combination of this decomposition with a weighted Osgood argument closes the desired estimate. Finally, the variational analysis requires a
different mechanism because coercivity
degenerates at the critical mass and the logarithmic Green charge breaks
mass-critical scaling.  The normalized bound state gives strict binding below
\(\mu_*\), while a mass-preserving Green charge variation and the logarithmic
scaling anomaly imply that the infimum of energy at the free critical mass is \(-\infty\).  Near the critical mass, we
establish a concentration-scale decomposition showing that the rescaled Green
charge vanishes, which yields rigidity toward free ground states.

The rest of the paper is organized as follows.  Section~\ref{sec:operator} introduces
the working space and basic estimates.
Section~\ref{sec:cubic-linear} develops the weighted linear estimates and
proves unconditional uniqueness.  Section~\ref{sec:cubic-flow} constructs the
strong flow and establishes global continuation below the critical mass.
Section~\ref{sec:cubic-ground-states} establishes the existence and orbital
stability of fixed-mass ground states and identifies the sharp critical threshold.
Finally, in Section~\ref{sec:critical-mass-rigidity}
we study near-critical rigidity toward free ground states.

\section{Preliminaries}
\label{sec:operator}

In this section, we collect the preliminaries from operator theory and
variational analysis that will be used throughout the paper.  We take the
\(L^2\)-inner product to be
\[
 \langle f,g\rangle:=\int_\R\overline{f(x)}g(x)\,\dd x,
\]
conjugate-linear in the first entry and linear in the second.  All Hilbert space inner products,
form pairings and dual pairings use the same convention.  If
\(X\hookrightarrow L^2\) is a form space, then
\(X^*\) denotes its anti-dual relative to the \(L^2\)-pivot.
For \(1\leq p<\infty\), we write
\[
 \norm f_p:=\left(\int_\R\abs{f(x)}^p\,\dd x\right)^{1/p},
 \qquad
 \norm f_\infty:=\mathop{\mathrm{ess\,sup}}_{x\in\R}\abs{f(x)}.
\]

Writing \(u=\phi+qG_{2\lambda_\alpha}\), the norm in
\eqref{eq:intro-form-norm} has the explicit form
\begin{equation}\label{eq:critical-form-norm}
 \norm u_{X_\alpha}^2
 =\norm{|D|^{1/2}\phi}_2^2
  +2\lambda_\alpha\norm\phi_2^2+\frac{\log2}{\pi}|q|^2.
\end{equation}
We have the continuous embeddings
\[
 X_\alpha\hookrightarrow L^2(\R)\hookrightarrow X_\alpha^*.
\]
The next lemma combines the case \(s=1\) of
\cite[Theorem~3.1]{MichelangeliScandone} with the free resolvent identity.

\begin{lemma}
\label{lem:operator-regimes}
The operator \(|D|_\alpha\) is bounded from below, and
for \(\lambda>0\), \(\lambda\neq\lambda_\alpha\) and \(f\in L^2(\R)\),
\begin{equation}\label{eq:critical-resolvent}
 (|D|_\alpha+\lambda)^{-1}f
 =(|D|+\lambda)^{-1}f
 +\frac{\langle G_\lambda,f\rangle}{\Phi_\alpha(-\lambda)}G_\lambda.
\end{equation}
Moreover, the following spectral properties hold.
\begin{enumerate}[\upshape(i)]
\item \(\sigma_{\mathrm{ess}}(|D|_\alpha)=[0,\infty)\).
\item \(\sigma(|D|_\alpha)\cap(-\infty,0)=\{-\lambda_\alpha\}\).
\item \(\ker(|D|_\alpha+\lambda_\alpha)
=\operatorname{span}\{G_{\lambda_\alpha}\}\).
\end{enumerate}
\end{lemma}

\begin{lemma}\label{lem:critical-green-bounds}
For every \(\lambda>0\), there is a constant \(C_\lambda>0\) such that
\begin{equation}\label{eq:critical-green-near}
 \abs{G_\lambda(x)-\frac1\pi\log\frac1{\lambda|x|}}
 \leq C_\lambda
 \qquad(0<\lambda|x|\leq1),
\end{equation}
and
\begin{equation}\label{eq:critical-green-Lp}
 \norm{G_\lambda}_p\leq C_\lambda p
 \qquad(2\leq p<\infty).
\end{equation}
\end{lemma}

\begin{proof}
Starting from \eqref{eq:intro-half-wave-green}, the Laplace transform and the
Poisson kernel give, for \(x\neq0\),
\begin{equation}\label{eq:critical-green-poisson}
 G_\lambda(x)
 =\frac1\pi\int_0^\infty
 e^{-\lambda t}\frac{t}{t^2+x^2}\,\dd t.
\end{equation}
Set \(a=\lambda|x|\), scale by
\(t=|x|s\) and split the integral at \(1\) and \(a^{-1}\).  On
\([1,a^{-1}]\), subtract \(s^{-1}\).  Since \(1-e^{-as}\leq as\), the two
resulting errors satisfy
\[
 \abs{\frac{s}{1+s^2}-\frac1s}\leq Cs^{-3}
 \quad(s\geq1),
 \qquad
 \int_1^{a^{-1}}\frac{1-e^{-as}}s\,\dd s\leq C.
\]
For the remaining pieces, use
\(s/(1+s^2)\leq\min\{s,s^{-1}\}\) and \(s^{-1}\leq a\) on
\([a^{-1},\infty)\).  Then
\[
 \int_0^1 e^{-as}\frac{s}{1+s^2}\,\dd s\leq C,
 \qquad
 \int_{a^{-1}}^\infty e^{-as}\frac{s}{1+s^2}\,\dd s
 \leq C.
\]
This proves \eqref{eq:critical-green-near}.  Moreover,
\eqref{eq:critical-green-poisson} gives
\[
 G_\lambda(x)\leq\frac1{\pi\lambda^2x^2}.
\]
Combining this far-field estimate with \eqref{eq:critical-green-near} yields
\[
 \norm{G_\lambda}_p^p
 \leq C_\lambda^p\bigl(1+\Gamma(p+1)\bigr),
\]
where \(\Gamma(s):=\int_0^\infty t^{s-1}e^{-t}\,\dd t\) for \(s>0\).
Since
\[
 \Gamma(p+1)^{1/p}\leq Cp,
\]
we obtain \eqref{eq:critical-green-Lp}.
\end{proof}

\begin{lemma}\label{lem:green-coordinates}
Let \(u\in X_\alpha\) and \(\lambda,\mu>0\).  Then the charge \(q\) is
independent of the Green parameter, which means that
\[
 u=\phi_\lambda+qG_\lambda=\phi_\mu+qG_\mu,
\]
and the Green functions satisfy
\[
 G_\lambda-G_\mu
 =(\mu-\lambda)(|D|+\lambda)^{-1}G_\mu\in H^{1/2}(\R).
\]
Moreover, for \(\rho>0\),
\[
 \D_\rho u=\D_\rho\phi_\lambda+\rho^{1/2}qG_{\rho\lambda}.
\]
\end{lemma}

\begin{proof}
The definition of \(G_\lambda\) gives
\[
 \widehat{G_\lambda-G_\mu}(\xi)
 =\frac{\mu-\lambda}{(|\xi|+\lambda)(|\xi|+\mu)}.
\]
This is the Fourier transform of
\((\mu-\lambda)(|D|+\lambda)^{-1}G_\mu\).  It is bounded near the origin
and is \(O(|\xi|^{-2})\) at infinity.  Hence
\(G_\lambda-G_\mu\in H^{1/2}(\R)\).

We can therefore rewrite the \(\lambda\)-decomposition as
\[
 u=\bigl[\phi_\lambda+q(G_\lambda-G_\mu)\bigr]+qG_\mu.
\]
The term in brackets belongs to \(H^{1/2}(\R)\).  This is the
\(\mu\)-decomposition of \(u\).  Its uniqueness shows that its charge is
again \(q\).

Finally, the Fourier representation gives
\(G_\lambda(\rho x)=G_{\rho\lambda}(x)\).  Thus
\[
 \D_\rho u
 =\D_\rho\phi_\lambda+\rho^{1/2}qG_\lambda(\rho x)
 =\D_\rho\phi_\lambda+\rho^{1/2}qG_{\rho\lambda}.
\]
\end{proof}

\begin{lemma}
\label{lem:critical-form-coordinates}
For every \(2\leq p<\infty\), the embedding
\(X_\alpha\hookrightarrow L^p(\R)\) is continuous and is compact on every
bounded interval.  Moreover, with the decomposition
\(u=\phi+qG_{2\lambda_\alpha}\),
\begin{equation}\label{eq:critical-regular-Lp}
 \norm\phi_p\leq C_\alpha\sqrt p\,\norm u_{X_\alpha}.
\end{equation}
\end{lemma}

\begin{proof}
Write \(u=\phi+qG_{2\lambda_\alpha}\).  The norm identity
\eqref{eq:critical-form-norm} controls \(\norm\phi_{H^{1/2}}+|q|\) by
\(\norm u_{X_\alpha}\).  Let \(p'=p/(p-1)\) and, for \(p>2\), set
\(r=2p/(p-2)\), so that \(1/p'=1/2+1/r\).  Hausdorff--Young and weighted
H\"older give
\[
 \norm\phi_p
 \leq C\norm{\widehat\phi}_{p'}
 \leq C\norm{(1+\abs\xi)^{1/2}\widehat\phi}_2
       \norm{(1+\abs\xi)^{-1/2}}_r.
\]
Moreover,
\[
 \int_\R(1+\abs\xi)^{-r/2}\,\dd\xi=p-2,
 \qquad
 (p-2)^{1/r}\leq\sqrt p.
\]
The case \(p=2\) follows directly from Plancherel's identity.  The coordinate
control therefore proves \eqref{eq:critical-regular-Lp}.  Together with
Lemma~\ref{lem:critical-green-bounds}, it also gives
\[
 \norm u_p
 \leq\norm\phi_p+|q|\norm{G_{2\lambda_\alpha}}_p
 \leq C_{\alpha,p}\norm u_{X_\alpha},
\]
which proves the continuous embedding.  For a bounded sequence,
Rellich compactness applies to the regular coordinates and finite
dimensionality to the charges.  Hence a subsequence converges strongly in
\(L^2(K)\).  For fixed \(p\geq2\), choose \(r>p\).  Interpolation between this
local convergence and the uniform \(L^r\)-bound from the continuous embedding
gives strong convergence in \(L^p(K)\).
\end{proof}

We use \(\alpha=\infty\) only to denote the free extension and set
\(|D|_\infty:=|D|\).  For \(v\in H^{1/2}(\R)\), define the free energy and its
fixed-mass level by
\[
 E_\infty(v)
 :=\frac12\norm{|D|^{1/2}v}_2^2-\frac14\norm v_4^4,
 \qquad
 I_\infty(\mu)
 :=\inf_{\substack{v\in H^{1/2}(\R)\\ M(v)=\mu}}E_\infty(v).
\]

\begin{lemma}
\label{lem:free-cubic-geometry}
The sharp critical inequality is
\begin{equation}\label{eq:intro-half-wave-GN}
 \norm v_4^4\leq\frac{2}{\mu_*}
 \norm{|D|^{1/2}v}_2^2\norm v_2^2,
 \qquad v\in H^{1/2}(\R).
\end{equation}
Moreover, the free energy levels satisfy
\begin{equation}\label{eq:cubic-free-levels}
 I_\infty(\mu)=0\quad(0<\mu\leq\mu_*),
 \qquad
 I_\infty(\mu)=-\infty\quad(\mu>\mu_*).
\end{equation}
The infimum is attained only at \(\mu=\mu_*\), precisely by the free
optimizers.
\end{lemma}

\begin{proof}
The existence of optimizers for the fractional Gagliardo--Nirenberg inequality
follows from \cite{BellazziniFrankVisciglia}.  For the cubic half-wave problem,
the corresponding sharp variational formulation is given in
\cite{KriegerLenzmannRaphael}.  Under the normalization
\eqref{eq:intro-free-ground-state}, the Pohozaev identities give the optimal
constant \(2/\mu_*\), while the optimizer classification follows from
\cite{FrankLenzmann}.  If \(M(v)=\mu\), then
\[
 E_\infty(v)
 \geq\frac12\left(1-\frac{\mu}{\mu_*}\right)
 \norm{|D|^{1/2}v}_2^2.
\]
Moreover, the mass-preserving dilation satisfies
\[
 E_\infty(\D_\rho v)=\rho E_\infty(v).
\]
These two facts give \(I_\infty(\mu)=0\) for
\(0<\mu\leq\mu_*\), and strict positivity in the first display excludes
attainment when \(0<\mu<\mu_*\).  If \(\mu>\mu_*\) and
\(a^2=\mu/\mu_*\), the Pohozaev identities for \(Q_*\) give
\[
 E_\infty(aQ_*)
 =\frac{\mu}{2}\left(1-\frac{\mu}{\mu_*}\right)<0.
\]
Letting \(\rho\to\infty\) yields \(I_\infty(\mu)=-\infty\).  At
\(\mu=\mu_*\), attainment and the optimizer statement follow from the
equality cases in \eqref{eq:intro-half-wave-GN}.
\end{proof}

The following coercivity estimate provides the uniform bound needed for
global continuation and will also be used in the constrained variational
analysis of Section~\ref{sec:cubic-ground-states}.

\begin{lemma}
\label{lem:cubic-half-wave-coercivity}
For every \(m\in(0,\mu_*)\), there are \(c_m,C_m>0\) such that
\begin{equation}\label{eq:cubic-half-wave-coercivity}
 E_\alpha(u)+C_m
 \geq c_m\norm u_{X_\alpha}^2,
 \qquad M(u)\leq m.
\end{equation}
\end{lemma}

\begin{proof}
We first prove the asymptotic estimate
\[
 \limsup_{n\to\infty}
 \frac{\norm{u_n}_4^4}{\norm{u_n}_{X_\alpha}^2}
 \leq\frac{2m_+}{\mu_*}
\]
whenever \(M(u_n)\leq m_+<\infty\) and
\(\norm{u_n}_{X_\alpha}\to\infty\).
Put \(\tau_n:=\mathfrak t_\alpha[u_n]\).  Since the masses are bounded,
\eqref{eq:critical-form-norm} shows that \(\tau_n\to\infty\) and
\(\tau_n/\norm{u_n}_{X_\alpha}^2\to1\).
For large \(n\), set
\[
 \Lambda_n:=\frac{\tau_n}{\sqrt{\log\tau_n}},
\]
so that \(\Lambda_n\to\infty\) and, in particular,
\(\Lambda_n>\lambda_\alpha\) for all sufficiently large \(n\).  This choice
simultaneously has \(\Lambda_n=o(\tau_n)\) and
\(\log\Lambda_n\sim\log\tau_n\): the mass shift is lower order while the
logarithmic coefficient still controls the charge.
Lemma~\ref{lem:green-coordinates} allows us to use the decomposition
\[
 u_n=\phi_n+q_nG_{\Lambda_n},\qquad
 \phi_n\in H^{1/2}(\R).
\]
Since
\[
 \Phi_\alpha(-\Lambda_n)
 =\frac1\pi\log\frac{\Lambda_n}{\lambda_\alpha}>0,
\]
the renormalized form identity gives
\begin{equation}\label{eq:cubic-half-wave-moving-decomposition}
 \tau_n+\Lambda_n\norm{u_n}_2^2
 =\norm{|D|^{1/2}\phi_n}_2^2+\Lambda_n\norm{\phi_n}_2^2
  +\Phi_\alpha(-\Lambda_n)|q_n|^2.
\end{equation}
The masses are bounded and
\[
 \Lambda_n\norm{u_n}_2^2
 =O\!\left(\frac{\tau_n}{\sqrt{\log\tau_n}}\right)=o(\tau_n).
\]
All terms on the right side of
\eqref{eq:cubic-half-wave-moving-decomposition} are nonnegative.  Consequently,
\begin{equation}\label{eq:cubic-half-wave-moving-bounds}
 \norm{|D|^{1/2}\phi_n}_2^2=O(\tau_n),
 \qquad
 |q_n|^2=O\left(\frac{\tau_n}{\log\tau_n}\right).
\end{equation}

The exact scaling \(G_\Lambda(x)=G_1(\Lambda x)\) yields
\[
 \norm{G_\Lambda}_2^2=\Lambda^{-1}\norm{G_1}_2^2,
 \qquad
 \norm{G_\Lambda}_4^4=\Lambda^{-1}\norm{G_1}_4^4.
\]
Estimates \eqref{eq:cubic-half-wave-moving-bounds} therefore imply
\begin{equation}\label{eq:cubic-half-wave-charge-smallness}
 \norm{q_nG_{\Lambda_n}}_2^2
 =O((\log\tau_n)^{-1/2})=o(1).
\end{equation}
\[
 \norm{q_nG_{\Lambda_n}}_4^4
 =O\left(\frac{\tau_n}{(\log\tau_n)^{3/2}}\right)=o(\tau_n).
\]
Writing \(g_n=q_nG_{\Lambda_n}\), the boundedness of
\(\norm{u_n}_2\) and \(\norm{g_n}_2=o(1)\) gives
\[
 \abs{\norm{\phi_n}_2^2-\norm{u_n}_2^2}
 \leq2\norm{u_n}_2\norm{g_n}_2+\norm{g_n}_2^2=o(1).
\]
Thus
\[
 \norm{\phi_n}_2^2=\norm{u_n}_2^2+o(1)\leq m_++o(1).
\]
By \eqref{eq:intro-half-wave-GN} and
\eqref{eq:cubic-half-wave-moving-decomposition},
\[
 \norm{\phi_n}_4^4
 \leq\frac{2\norm{\phi_n}_2^2}{\mu_*}
       \norm{|D|^{1/2}\phi_n}_2^2=O(\tau_n).
\]
The elementary quartic difference estimate and H\"older give
\[
 \abs{\norm{\phi_n+g_n}_4^4-\norm{\phi_n}_4^4}
 \leq C\sum_{j=1}^4
 \norm{\phi_n}_4^{4-j}\norm{g_n}_4^j=o(\tau_n),
\]
 because \(\norm{\phi_n}_4^4=O(\tau_n)\) and
\(\norm{g_n}_4^4=o(\tau_n)\).
Hence
\[
 \norm{u_n}_4^4=\norm{\phi_n}_4^4+o(\tau_n).
\]
Moreover, since \(\Lambda_n=o(\tau_n)\),
\eqref{eq:cubic-half-wave-moving-decomposition} gives
\[
 \norm{|D|^{1/2}\phi_n}_2^2\leq\tau_n+o(\tau_n).
\]
Combining the last three estimates proves
\[
 \limsup_{n\to\infty}\frac{\norm{u_n}_4^4}{\tau_n}
 \leq\frac{2m_+}{\mu_*},
\]
and hence the asserted asymptotic estimate.

To obtain cubic coercivity, consider any sequence with
\(M(u_n)\leq m<\mu_*\) and \(\norm{u_n}_{X_\alpha}\to\infty\), and apply the
asymptotic estimate with \(m_+=m\).  It gives
\[
 \liminf_{n\to\infty}
 \frac{E_\alpha(u_n)}
      {\norm{u_n}_{X_\alpha}^2}
 \geq\frac12-\frac{m_+}{2\mu_*}>0.
\]
Choose
\[
 0<c_m<\frac12-\frac{m}{2\mu_*}.
\]
If no radius gave
\(E_\alpha(u)\geq c_m\norm u_{X_\alpha}^2\) for all
\(M(u)\leq m\) outside that radius, a sequence of counterexamples with form
norm tending to infinity would contradict the preceding liminf estimate.
The energy is therefore coercive outside a fixed ball in the form norm.  Inside
that ball, the continuous embedding
\(X_\alpha\hookrightarrow L^4\) and
\[
 \abs{\mathfrak t_\alpha[u]}
 \leq \norm u_{X_\alpha}^2+2\lambda_\alpha\norm u_2^2
\]
give the explicit uniform bound
\[
 \abs{E_\alpha(u)}\leq C_{\alpha,m,R}
 \qquad
 \bigl(\norm u_{X_\alpha}\leq R,\ M(u)\leq m\bigr).
\]
Increasing the additive constant to cover this ball proves
\eqref{eq:cubic-half-wave-coercivity}.
\end{proof}

\section{Unconditional uniqueness in the energy space}
\label{sec:cubic-linear}

In this section we establish unconditional uniqueness in the natural energy
space. We recall the cubic half-wave equation with a point interaction
\begin{equation}\label{eq:cubic-evolution-recalled}
 i\partial_t u=|D|_\alpha u-|u|^{2}u,\qquad (t,x)\in\R\times\R.
\end{equation}
We begin with the notion of solution used throughout the argument.

\begin{definition}[Weak energy solution]
\label{def:cubic-half-wave-weak-solution}
Let \(J\subset\R\) be an interval.  A weak energy solution of
\eqref{eq:cubic-evolution-recalled} on \(J\) is a function
\begin{equation}\label{eq:cubic-half-wave-weak-class}
 u\in C_w(J;X_\alpha)\cap C(J;L^2)
 \cap W^{1,\infty}_{\mathrm{loc}}(J;X_\alpha^*)
 \cap L^\infty_{\mathrm{loc}}(J;X_\alpha)
\end{equation}
such that, for almost every \(t\in J\) and every \(v\in X_\alpha\),
\[
 \langle i\partial_tu(t),v\rangle_{X_\alpha^*,X_\alpha}
 =\mathfrak t_\alpha(u(t),v)
  -\langle |u(t)|^2u(t),v\rangle_{X_\alpha^*,X_\alpha}.
\]
\end{definition}

\begin{definition}[Unconditional uniqueness]
\label{def:cubic-half-wave-unconditional-uniqueness}
Let \(J\subset\R\) be an interval,
and suppose that \(u\) and \(v\) are weak energy solutions on \(J\). We say that weak energy solutions are unconditionally
unique if
\[
 u(t_0)=v(t_0)\quad\text{in }L^2(\R)
\]
for some \(t_0\in J\) implies \(u=v\) on \(J\).
\end{definition}

The main result of this section is the following proposition.

\begin{proposition}
\label{prop:cubic-half-wave-weak-unique}
Weak energy solutions of \eqref{eq:evolution} are unconditionally unique.
\end{proposition}

The rest of this section is devoted to the proof of
Proposition~\ref{prop:cubic-half-wave-weak-unique}.  We first establish a
localized Kato smoothing estimate on the continuous spectrum.  A Green
function decomposition then reduces uniqueness to a weighted Osgood argument.

\subsection{Localized Kato smoothing}

For \(z\in\C\setminus[0,\infty)\), define \(\mathcal G_z\in L^2(\R)\) by
\begin{equation}\label{eq:cubic-complex-green}
 \widehat{\mathcal G_z}(\xi)=\frac1{\abs\xi-z}.
\end{equation}
For \(x\neq0\), its physical space representative is
\[
 \mathcal G_z(x)=\frac1\pi\int_0^\infty
 \frac{\cos(\xi x)}{\xi-z}\,\dd\xi,
\]
where the integral is understood as an improper oscillatory integral.  In
particular, \(\mathcal G_{-\lambda}=G_\lambda\) for \(\lambda>0\).
If \(z=re^{i\theta}\), where \(r>0\) and \(0<\theta<2\pi\), define
\[
\Phi_\alpha(z)
 :=\frac1\pi\log\frac r{\lambda_\alpha}
   +\frac{i}{\pi}(\theta-\pi).
\]
Let \(P_{\mathrm c}\) denote the orthogonal projection onto the continuous
spectral subspace of \(|D|_\alpha\).  The normalized eigenfunction and the
projection are given by
\begin{equation}\label{eq:cubic-continuous-projection}
 h_\alpha
 :=\frac{G_{\lambda_\alpha}}{\norm{G_{\lambda_\alpha}}_2},
 \qquad
 P_{\mathrm c}f:=f-\langle h_\alpha,f\rangle h_\alpha.
\end{equation}
The corresponding projected resolvent kernel is given by the following
lemma.

\begin{lemma}\label{lem:cubic-point-resolvent-kernels}
Let
\(z\in\C\setminus([0,\infty)\cup\{-\lambda_\alpha\})\).  Then, for almost every
\(x,y\in\R\), the kernel of the projected resolvent is
\begin{equation}\label{eq:cubic-Pc-kernel}
 \bigl[P_{\mathrm c}(|D|_\alpha-z)^{-1}P_{\mathrm c}\bigr](x,y)
={}\mathcal G_z(x-y)
 +\frac{\mathcal G_z(x)\mathcal G_z(y)}{\Phi_\alpha(z)}-\frac{\pi\lambda_\alpha
 G_{\lambda_\alpha}(x)G_{\lambda_\alpha}(y)}{-\lambda_\alpha-z}.
\end{equation}
\end{lemma}

\begin{proof}
The map \(z\mapsto\mathcal G_z\) is \(L^2\)-analytic.  For
\(z=-\lambda\neq-\lambda_\alpha\), formula
\eqref{eq:critical-resolvent} gives
\begin{equation}\label{eq:cubic-complex-Krein-operator}
 (|D|_\alpha-z)^{-1}f
 =(|D|-z)^{-1}f
 +\frac{\langle\mathcal G_{\overline z},f\rangle}
 {\Phi_\alpha(z)}\mathcal G_z,
 \qquad f\in L^2(\R).
\end{equation}
Both sides admit meromorphic continuations to
\(\C\setminus[0,\infty)\).  The operator-valued identity theorem therefore
proves \eqref{eq:cubic-complex-Krein-operator} throughout the resolvent set.
Since \(\Phi_\alpha\) has a simple zero at \(-\lambda_\alpha\), the rank-one
term has a simple pole there.

The inner product is linear in its second argument, and
\[
 \overline{\mathcal G_{\overline z}(y)}=\mathcal G_z(y).
\]
Hence \eqref{eq:cubic-complex-Krein-operator} has the kernel
\begin{equation}\label{eq:cubic-complex-Krein}
 (|D|_\alpha-z)^{-1}(x,y)
 =\mathcal G_z(x-y)
 +\frac{\mathcal G_z(x)\mathcal G_z(y)}{\Phi_\alpha(z)}.
\end{equation}
Plancherel's identity gives
\[
 \norm{G_{\lambda_\alpha}}_2^2
 =\frac1{2\pi}\int_\R\frac{\dd\xi}
 {(\abs\xi+\lambda_\alpha)^2}
 =\frac1{\pi\lambda_\alpha}.
\]
Since \(I-P_{\mathrm c}\) is the orthogonal projection onto
\(\operatorname{span}\{G_{\lambda_\alpha}\}\), it has kernel
\begin{equation}\label{eq:cubic-bound-kernel}
 (I-P_{\mathrm c})(x,y)
 =\pi\lambda_\alpha
 G_{\lambda_\alpha}(x)G_{\lambda_\alpha}(y).
\end{equation}
Finally, the spectral decomposition at the eigenvalue
\(-\lambda_\alpha\) gives
\[
 \bigl[P_{\mathrm c}(|D|_\alpha-z)^{-1}P_{\mathrm c}\bigr](x,y)
 =(|D|_\alpha-z)^{-1}(x,y)
 -\frac{(I-P_{\mathrm c})(x,y)}{-\lambda_\alpha-z}.
\]
Combining this identity with \eqref{eq:cubic-complex-Krein} and
\eqref{eq:cubic-bound-kernel} proves \eqref{eq:cubic-Pc-kernel}.
\end{proof}

To isolate the threshold cancellation used below, set
\(\widetilde{\mathcal G}_z(t):=\mathcal G_z(t)+\Phi_\alpha(z)-\alpha\).
A direct substitution into \eqref{eq:cubic-complex-Krein} gives
\begin{equation}\label{eq:cubic-Krein-cancellation}
 (|D|_\alpha-z)^{-1}(x,y)
={}\widetilde{\mathcal G}_z(x-y)-\widetilde{\mathcal G}_z(x)
 -\widetilde{\mathcal G}_z(y)-\alpha+\frac{(\alpha+\widetilde{\mathcal G}_z(x))
 (\alpha+\widetilde{\mathcal G}_z(y))}
 {\Phi_\alpha(z)}.
\end{equation}

We first record the free kernel estimates used below.  Set
\[
 \ell(x):=1+\log_+\frac1{\abs x},\qquad
 \log_+a:=\max\{\log a,0\}.
\]
\begin{lemma}
\label{lem:cubic-free-boundary-kernels}
For \(s>0\) and \(t\neq0\), the limits
\[
 \mathcal G_{s\pm i0}(t)
 :=\lim_{\varepsilon\to0^+}\mathcal G_{s\pm i\varepsilon}(t)
\]
exist and satisfy
\begin{equation}\label{eq:cubic-free-boundary-bounds}
 \abs{\mathcal G_{-s}(t)}+\abs{\mathcal G_{s\pm i0}(t)}
 \leq C\left(1+\log_+\frac1{s\abs t}\right).
\end{equation}
\end{lemma}

\begin{proof}
Using this representation and changing variables by \(\eta=\xi\abs t\),
\[
 \pi\mathcal G_{-s}(t)
 =\int_0^\infty\frac{\cos\eta}{\eta+s\abs t}\,\dd\eta.
\]
The standard identity
\[
 \int_0^\infty\frac{\cos\eta}{\eta+a}\,\dd\eta
 =-\cos a\,\operatorname{Ci}(a)
  +\sin a\left(\frac\pi2-\operatorname{Si}(a)\right)
\]
holds for \(a>0\), where
\[
 \operatorname{Ci}(a):=-\int_a^\infty\frac{\cos\eta}{\eta}\,\dd\eta
 \quad\text{and}\quad
 \operatorname{Si}(a):=\int_0^a\frac{\sin\eta}{\eta}\,\dd\eta.
\]
The expansions of \(\operatorname{Ci}\) and \(\operatorname{Si}\) at
zero, followed by one integration by parts for \(a\geq1\), prove
the negative part of \eqref{eq:cubic-free-boundary-bounds}.

On the positive boundary, the Sokhotski--Plemelj formula gives
\[
 \pi\mathcal G_{s\pm i0}(t)
 =\operatorname{p.v.}\int_0^\infty
 \frac{\cos\eta}{\eta-s\abs t}\,\dd\eta
 \mathbin{\pm}i\pi\cos(s\abs t).
\]
Using
\[
 \operatorname{p.v.}\int_0^\infty\frac{\cos\eta}{\eta-a}\,\dd\eta
 =-\cos a\,\operatorname{Ci}(a)
   -\sin a\left(\operatorname{Si}(a)+\frac\pi2\right)
\]
in the same way proves the positive part of
\eqref{eq:cubic-free-boundary-bounds}.
\end{proof}

For \(s>0\), the definition of \(\Phi_\alpha\) gives
\[
 \Phi_\alpha(s\pm i0)
 :=\lim_{\varepsilon\to0^+}\Phi_\alpha(s\pm i\varepsilon)
 =\frac1\pi\log\frac{s}{\lambda_\alpha}\mathbin{\mp}i.
\]
Lemma~\ref{lem:cubic-free-boundary-kernels} then shows that the limits
\[
 \widetilde{\mathcal G}_{s\pm i0}(t)
 :=\lim_{\varepsilon\to0^+}
 \widetilde{\mathcal G}_{s\pm i\varepsilon}(t)
 =\mathcal G_{s\pm i0}(t)+\Phi_\alpha(s\pm i0)-\alpha
\]
exist for \(t\neq0\).  For \(0<\abs t\leq1\) and
\(0<s\abs t\leq1\), the small argument expansions above and
\(\lambda_\alpha=e^{-\pi\alpha}\) give
\begin{equation}\label{eq:cubic-renormalized-boundary-bounds}
 \abs{\widetilde{\mathcal G}_{-s}(t)}
 +\abs{\widetilde{\mathcal G}_{s\pm i0}(t)}
 \leq C\ell(t).
\end{equation}
For every \(t\neq0\), they also give
\begin{equation}\label{eq:cubic-renormalized-threshold}
 \lim_{s\to0^+}\widetilde{\mathcal G}_{s\pm i0}(t)
 =-\frac1\pi\bigl(\gamma+\log\abs t\bigr).
\end{equation}
Here \(\gamma\) is the Euler--Mascheroni constant, given by
\(\gamma=-\int_0^\infty e^{-r}\log r\,\dd r\).

We now pass to the projected resolvent kernel.  For
\(\sigma\in\R\), let
\[
 \mathcal R_{\alpha,\mathrm c}^{\pm}(\sigma;x,y)
 :=\bigl[P_{\mathrm c}(|D|_\alpha-(\sigma\pm i0))^{-1}
 P_{\mathrm c}\bigr](x,y)
\]
denote the projected resolvent kernel.  For \(\sigma>0\), it is the
pointwise boundary value obtained from \eqref{eq:cubic-Pc-kernel}.  For
\(\sigma<0\), \(\sigma\neq-\lambda_\alpha\), it is the ordinary resolvent
kernel.  At \(\sigma=0\), it denotes the corresponding threshold limit.
After projection, the pole at \(\sigma=-\lambda_\alpha\) is removable, and
the value there is defined by analytic continuation.

\begin{lemma}
\label{lem:cubic-boundary-kernel}
There are \(r_*>0\) and \(C_\alpha>0\) such that, for almost every \(x,y\)
with \(\abs x,\abs y<2r_*\),
\begin{equation}\label{eq:cubic-kernel-majorant}
 \sup_{\sigma\in\R}
 \abs{\mathcal R_{\alpha,\mathrm c}^{\pm}(\sigma;x,y)}
 \leq C_\alpha\bigl(
  1+\ell(x-y)+\ell(x)+\ell(y)+\ell(x)\ell(y)\bigr).
\end{equation}
\end{lemma}

\begin{proof}
We insert the bounds of Lemma~\ref{lem:cubic-free-boundary-kernels} into the
continuous resolvent formula of
Lemma~\ref{lem:cubic-point-resolvent-kernels}.  Choose \(r_*<1/8\) so that
\(4r_*\lambda_\alpha<1/2\).  On the positive axis, the explicit formula
for \(\Phi_\alpha(s\pm i0)\) gives
\[
 \abs{\Phi_\alpha(s\pm i0)}
 \geq c\left(1+\abs{\log\frac{s}{\lambda_\alpha}}\right).
\]
If \(s\leq(4r_*)^{-1}\), use
\eqref{eq:cubic-renormalized-boundary-bounds} in
\eqref{eq:cubic-Krein-cancellation}.  If \(s>(4r_*)^{-1}\), use
\eqref{eq:cubic-complex-Krein}, \eqref{eq:cubic-free-boundary-bounds} and
\[
 \log_+\frac1{s\abs t}
 \leq\log_+\frac{4r_*}{\abs t}\leq C+\ell(t),
 \qquad \abs t<4r_*.
\]
 The finite-rank term in \eqref{eq:cubic-Pc-kernel} is bounded by
 \(C_\alpha\ell(x)\ell(y)\).  This proves
 \eqref{eq:cubic-kernel-majorant} on the positive axis.  At the threshold,
 \eqref{eq:cubic-renormalized-threshold} and
 \(1/\Phi_\alpha(s\pm i0)\to0\) give the same bound.

On the negative axis write \(z=-s\) and assume
\(\abs{s-\lambda_\alpha}\geq\lambda_\alpha/2\).  Then
\[
 \Phi_\alpha(-s)=\frac1\pi\log\frac{s}{\lambda_\alpha}
\]
is separated from zero.  For \(s\leq(4r_*)^{-1}\), use
\eqref{eq:cubic-Krein-cancellation} and
\eqref{eq:cubic-renormalized-boundary-bounds}.  For larger \(s\), use
 \eqref{eq:cubic-complex-Krein} and
 \eqref{eq:cubic-free-boundary-bounds}.  The
 finite-rank term in \eqref{eq:cubic-Pc-kernel} satisfies the same majorant
 because \(\abs{s-\lambda_\alpha}\geq\lambda_\alpha/2\).

 It remains to consider \(z=-s\) with
 \(\abs{s-\lambda_\alpha}<\lambda_\alpha/2\).  At the operator level,
 the spectral theorem gives
 \[
  P_{\mathrm c}(|D|_\alpha-z)^{-1}P_{\mathrm c}
  =\int_{[0,\infty)}\frac{1}{\lambda-z}\,
    \dd E_{|D|_\alpha}(\lambda),
 \]
 which is analytic through \(z=-\lambda_\alpha\).  We now identify a
 spatially uniform kernel majorant for this analytic value.  Since
 \[
  \Phi_\alpha(-\lambda_\alpha)=0,\qquad
  \left.\frac{\dd}{\dd s}\Phi_\alpha(-s)
  \right|_{s=\lambda_\alpha}
  =\frac1{\pi\lambda_\alpha},
 \]
 put \(\delta=s-\lambda_\alpha\),
 \(\mathcal K_s(x,y)=G_s(x)G_s(y)\) and
 \(a_\alpha=(\pi\lambda_\alpha)^{-1}\).  Differentiation under the absolutely
 convergent integral gives, for \(j=1,2\),
 \[
  \partial_s^jG_s(t)
  =\frac{(-1)^j j!}{\pi}\int_0^\infty
  \frac{\cos(\xi t)}{(\xi+s)^{j+1}}\,\dd\xi.
 \]
 Hence these derivatives are uniformly bounded on the pole interval.  The
 positive representation \eqref{eq:critical-green-poisson}, split at
 \(\tau=\abs t\) and \(\tau=1\), also gives
 \(\abs{G_s(t)}\leq C_\alpha\ell(t)\) uniformly on this interval.
 Taylor's formula with integral remainder now gives
 \[
  \Phi_\alpha(-s)
  =a_\alpha\delta+O_\alpha(\delta^2).
 \]
 \begin{align*}
  \mathcal K_s
  =\mathcal K_{\lambda_\alpha}
    +\delta\left.\partial_\eta \mathcal K_\eta\right|_{\eta=\lambda_\alpha}+O_\alpha\!\left(
  \delta^2[1+\ell(x)+\ell(y)+\ell(x)\ell(y)]\right).
 \end{align*}
 Indeed,
 \[
  \partial_s^2\mathcal K_s(x,y)
  =\partial_s^2G_s(x)G_s(y)
   +2\partial_s G_s(x)\partial_s G_s(y)
   +G_s(x)\partial_s^2G_s(y),
 \]
 and the preceding bounds imply
 \[
  \abs{\partial_s^2\mathcal K_s(x,y)}
  \leq C_\alpha\bigl(1+\ell(x)+\ell(y)+\ell(x)\ell(y)\bigr).
 \]
 Moreover, \eqref{eq:cubic-bound-kernel} gives
 \[
  \frac{\mathcal K_{\lambda_\alpha}(x,y)}{a_\alpha}
  =(I-P_{\mathrm c})(x,y).
 \]
 Dividing the Taylor expansions and cancelling their principal parts yields
 \[
  \left|
  \frac{G_s(x)G_s(y)}{\Phi_\alpha(-s)}
  -\frac{(I-P_{\mathrm c})(x,y)}{s-\lambda_\alpha}
  \right|
  \leq C_\alpha\bigl(
  1+\ell(x)+\ell(y)+\ell(x)\ell(y)\bigr).
 \]
 Adding the free kernel completes the proof of
 \eqref{eq:cubic-kernel-majorant}, including the analytic value at the removed
 pole.
\end{proof}

By \(G_{2\lambda_\alpha}(x)\to+\infty\) as
\(x\to0\), there exists \(0<r_0<r_*\) such that
\begin{equation}\label{eq:cubic-green-lower}
 G_{2\lambda_\alpha}(x)\geq1
 \qquad(0<\abs x<2r_0).
\end{equation}
We now define the weight used in the projected resolvent estimate.  Choose
\(\chi\in C_c^\infty(\R)\), \(0\leq\chi\leq1\), with \(\chi=1\) on
\([-1,1]\) and \(\operatorname{supp}\chi\subset[-2,2]\).  Then, for
\(0<r<r_0\), set
\begin{equation}\label{eq:cubic-resolvent-weight}
 W_r(x):=\chi(x/r)^{1/2}
 \abs{G_{2\lambda_\alpha}(x)}^2.
\end{equation}
Although \(W_r\) is unbounded at the origin, it belongs to \(L^p(\R)\) for
every \(1\leq p<\infty\) and satisfies
\begin{equation}\label{eq:cubic-weight-small}
 \norm{W_r}_2^2
 \lesssim r(1+\abs{\log r})^4\to0.
\end{equation}

Since \(W_r\) is unbounded, the weighted resolvent is initially understood as
a sesquilinear form.  Fix \(0<r<r_0\) and \(n\in\mathbb N\), and set
\(W_{r,n}:=\min\{W_r,n\}\).  Then
\(W_{r,n}\mathcal R_{\alpha,\mathrm c}^{\pm}(\sigma)W_{r,n}\) denotes the
integral operator with the corresponding weighted kernel.
The following lemma gives the truncated weighted resolvent estimate needed
below.

\begin{lemma}
\label{lem:cubic-weighted-boundary}
There are \(N\geq0\) and \(C_\alpha>0\), independent of \(r\) and \(n\),
such that, for every \(g\in L^2(\R)\),
\begin{equation}\label{eq:cubic-truncated-resolvent}
 \sup_{z\in\C\setminus\R}
 \norm{W_{r,n}P_{\mathrm c}(|D|_\alpha-z)^{-1}
 P_{\mathrm c}W_{r,n}g}_2
 \leq C_\alpha r(1+\abs{\log r})^N\norm g_2.
\end{equation}
\end{lemma}

\begin{proof}
We first identify the boundary values of the truncated resolvent.  Fix
\(\sigma>0\), set
\(z=\sigma\pm i\varepsilon\), and take
\(0<\varepsilon<\sigma/2\).  On \((0,2\sigma)\), write
\begin{align*}
 \int_0^{2\sigma}\frac{\cos(\xi t)}{\xi-z}\,\dd\xi
 ={}&\cos(\sigma t)
 \int_0^{2\sigma}\frac{\dd\xi}{\xi-z}+\int_0^{2\sigma}
 \frac{\cos(\xi t)-\cos(\sigma t)}{\xi-z}\,\dd\xi.
\end{align*}
After the change of variables \(s=\xi-\sigma\), the first integral is
\(\pm2i\arctan(\sigma/\varepsilon)\), hence is uniformly bounded.  The last
integrand has absolute value at most \(\abs t\), because
\(\abs{\xi-z}\geq\abs{\xi-\sigma}\).  On \((2\sigma,\infty)\), split at
\(P=\max\{2\sigma,\abs t^{-1}\}\).  Absolute integration up to \(P\),
followed by one integration by parts on \((P,\infty)\), gives
\begin{equation}\label{eq:cubic-complex-local-majorant}
 \sup_{0<\varepsilon<\sigma/2}
 \abs{\mathcal G_{\sigma\pm i\varepsilon}(t)}
 \leq C_{\sigma,r_0}\ell(t),
 \qquad 0<\abs t<4r_0.
\end{equation}
Indeed, the part before \(P\) is bounded by
\(C(1+\log_+(1/(\sigma\abs t)))\), whereas the integration-by-parts
remainder is bounded by \(C/(\abs t P)\leq C\).

For \(z\) in this nontangential neighborhood,
\(\Phi_\alpha(z)\) stays separated from zero, since
\(\abs{\Phi_\alpha(\sigma\pm i0)}\geq1\).  Thus
\eqref{eq:cubic-complex-Krein},
\eqref{eq:cubic-Pc-kernel} and
\eqref{eq:cubic-complex-local-majorant} imply, uniformly in
\(\varepsilon\),
\[
 \abs{\bigl[P_{\mathrm c}(|D|_\alpha-z)^{-1}P_{\mathrm c}\bigr](x,y)}
 \leq C_{\alpha,\sigma,r_0}
 \bigl(1+\ell(x-y)+\ell(x)\ell(y)\bigr)
\]
for \(\abs x,\abs y<2r_0\).  The Sokhotski--Plemelj formula applied to the
free oscillatory integral gives pointwise convergence to the boundary
kernels \(\mathcal R_{\alpha,\mathrm c}^{\pm}(\sigma;x,y)\) for almost every
\((x,y)\).
Multiplication by \(W_{r,n}(x)W_{r,n}(y)\), with
\[
 W_{r,n}(x)\leq W_r(x)
 \lesssim\mathbf1_{\{\abs x<2r\}}\ell(x)^2,
\]
 produces a square-integrable majorant on \(\R^2\).  Dominated convergence in
 \(L^2(\R^2)\), together with the Hilbert--Schmidt kernel estimate, proves,
 for every \(g\in L^2(\R)\),
\[
 \lim_{\varepsilon\to0^+}
 \norm{W_{r,n}P_{\mathrm c}
 (|D|_\alpha-(\sigma\pm i\varepsilon))^{-1}P_{\mathrm c}W_{r,n}
 g-W_{r,n}\mathcal R_{\alpha,\mathrm c}^{\pm}(\sigma)W_{r,n}g}_2
 =0.
\]
On the negative axis the same identification follows from ordinary resolvent
convergence away from \(-\lambda_\alpha\).  At \(-\lambda_\alpha\), the
eigenprojection in \(P_{\mathrm c}\) removes the simple pole, and
Lemma~\ref{lem:cubic-boundary-kernel} identifies the resulting analytic
continuation.  Thus the integral operators with kernels
\(W_{r,n}(x)\mathcal R_{\alpha,\mathrm c}^{\pm}(\sigma;x,y)W_{r,n}(y)\)
are the almost everywhere boundary values of the truncated weighted
resolvent on the real axis.

We now estimate these boundary operators.  For fixed nonnegative integers
\(a,b,c\), scaling \(x=rX\), \(y=rY\) and local integrability of finite
powers of \(\log_+(1/|t|)\) give
\begin{equation}\label{eq:cubic-log-product-integral}
 \int_{\abs x,\abs y<2r}
 \ell(x)^a\ell(y)^b\ell(x-y)^c\,\dd x\,\dd y
 \leq C_{a,b,c}r^2(1+\abs{\log r})^{a+b+c}.
\end{equation}
Indeed, after scaling, each logarithmic factor is bounded by
\(1+|\log r|+\log_+(1/|X|)\), with analogous bounds for \(Y\) and \(X-Y\),
and H\"older's inequality controls the resulting integrals.
Lemma~\ref{lem:cubic-boundary-kernel},
\eqref{eq:cubic-log-product-integral} and
\[
 W_{r,n}(x)\leq W_r(x)
 \lesssim\mathbf1_{\{\abs x<2r\}}\ell(x)^2
\]
therefore imply
\begin{equation}\label{eq:cubic-boundary-operator}
 \sup_{\sigma\in\R}
 \norm{W_{r,n}\mathcal R_{\alpha,\mathrm c}^{\pm}(\sigma)
 W_{r,n}g}_2
 \leq C_\alpha r(1+\abs{\log r})^N\norm g_2,
 \qquad g\in L^2(\R),
\end{equation}
uniformly in \(n\).

It remains to propagate this bound into the two half-planes.  Let
\(E_\alpha\) be the spectral resolution of \(|D|_\alpha\).  Fix unit vectors
\(f,g\in L^2(\R)\), and define the finite complex measure
\[
 \mu_{f,g}(\Omega)
 :=\left\langle P_{\mathrm c}W_{r,n}f,
 E_\alpha(\Omega)P_{\mathrm c}W_{r,n}g\right\rangle.
\]
The spectral theorem gives
\[
 F_{f,g}(z)
 :=\left\langle f,W_{r,n}P_{\mathrm c}(|D|_\alpha-z)^{-1}
 P_{\mathrm c}W_{r,n}g\right\rangle
 =\int_\R\frac{1}{\lambda-z}\,\dd\mu_{f,g}(\lambda).
\]
For every finite Borel partition \(\{\Omega_j\}\), orthogonality of the
spectral projections and Cauchy--Schwarz give
\[
 \sum_j\abs{\mu_{f,g}(\Omega_j)}
 \leq\norm{P_{\mathrm c}W_{r,n}f}_2
       \norm{P_{\mathrm c}W_{r,n}g}_2.
\]
Hence \(\mu_{f,g}\) has finite total variation.

We treat the upper half-plane first.  The Cayley map
\(\chi(z)=(z-i)/(z+i)\) sends it to \(\mathbb D\), with
\(\chi^{-1}(\zeta)=i(1+\zeta)/(1-\zeta)\).  Let \(\nu_{f,g}\) be the
pushforward under \(\lambda\mapsto(\lambda+i)/(\lambda-i)\) of the measure
\((\lambda-i)^{-1}\mu_{f,g}\).  Then
\[
 F_{f,g}(\chi^{-1}(\zeta))
 =(1-\zeta)\int_{\partial\mathbb D}
 \frac{1}{1-\eta\zeta}\,\dd\nu_{f,g}(\eta).
\]
Moreover,
\[
 \norm{\nu_{f,g}}_{\mathrm{TV}}
 \leq\norm{\mu_{f,g}}_{\mathrm{TV}}<\infty.
\]
Thus \cite[Theorem~3.5]{CimaMathesonRoss} implies that the disk Cauchy
transform, and hence \(F_{f,g}\circ\chi^{-1}\), belongs to \(H^q(\mathbb D)\)
for every \(0<q<1\).  In particular, it lies in the Smirnov class and has
nontangential boundary values almost everywhere.  For almost every
\(\sigma\in\R\), dominated convergence in the truncated kernel formula
identifies these values with the scalar matrix elements of
\(\mathcal R_{\alpha,\mathrm c}^{+}(\sigma)\).  The exceptional boundary point
corresponding to infinity is irrelevant.  Therefore
\eqref{eq:cubic-boundary-operator} and the Smirnov maximum principle
\cite[Chapter~II, Section~5]{Garnett} give the same bound throughout the upper
half-plane.  The conjugate Cayley map gives the lower half-plane estimate.
Taking the supremum over unit \(f\) and then scaling \(g\) proves
\eqref{eq:cubic-truncated-resolvent}.
\end{proof}

We can now pass from the projected boundary kernels to localized Kato
smoothing.

\begin{proposition}
\label{prop:cubic-Kato}
For every fixed \(\alpha\in\R\) and \(0<r<r_0\), there is
\(C_{\alpha,r}>0\) such that, for \(f\in L^2(\R)\) and
\(\Psi\in L^2(\R_t\times\R_x)\),
\begin{equation}\label{eq:cubic-hom-smooth}
 \norm{W_re^{-it|D|_\alpha}P_{\mathrm c}f}_{L^2(\R_t\times\R_x)}
 \leq C_{\alpha,r}\norm f_2.
\end{equation}
Moreover, interpreting the integral by its bounded Kato extension, we have
\begin{equation}\label{eq:cubic-ret-smooth}
 \norm{W_r\int_0^t e^{-i(t-\tau)|D|_\alpha}
 P_{\mathrm c}W_r\Psi(\tau)\,\dd\tau}_{L^2(\R_t\times\R_x)}
 \leq C_{\alpha,r}\norm \Psi_{L^2(\R_t\times\R_x)}.
\end{equation}
\end{proposition}

\begin{proof}
Lemma~\ref{lem:cubic-weighted-boundary} gives the required resolvent bound
uniformly in the truncation parameter.  We first remove this truncation and
then verify the closed-operator hypotheses of Kato's theory.
For \(f,g\in \mathrm{Dom}(W_r)\), which is dense in \(L^2\), dominated convergence
gives \(W_{r,n}f\to W_rf\) and \(W_{r,n}g\to W_rg\) in \(L^2\).
Consequently,
\[
\begin{split}
 \langle f,W_{r,n}P_{\mathrm c}(|D|_\alpha-z)^{-1}
 P_{\mathrm c}W_{r,n}g\rangle\to
 \langle P_{\mathrm c}W_rf,(|D|_\alpha-z)^{-1}
 P_{\mathrm c}W_rg\rangle.
\end{split}
\]
The uniform truncated bound shows that the limiting form is bounded by
\(C_\alpha r(1+\abs{\log r})^N\norm f_2\norm g_2\).
Since \(\mathrm{Dom}(W_r)\) is dense, it is represented by a unique bounded
operator \(T_r(z)\) on \(L^2\), characterized by
\[
 \langle f,T_r(z)g\rangle
 =\langle P_{\mathrm c}W_rf,(|D|_\alpha-z)^{-1}
 P_{\mathrm c}W_rg\rangle,
 \qquad f,g\in\mathrm{Dom}(W_r).
\]
Thus there is \(\varepsilon_r\leq C_\alpha r(1+\abs{\log r})^N\) such that
\begin{equation}\label{eq:cubic-epsilon}
 \sup_{z\in\C\setminus\R}
 \norm{T_r(z)g}_2
 \leq\varepsilon_r\norm g_2.
\end{equation}
We now connect this bounded form extension with the closed-operator
formulation of Kato smoothing.

The logarithmic behavior at zero gives the explicit majorant
\[
 \abs{W_rh_\alpha}\leq C_{\alpha,r}
 \mathbf1_{\{\abs x<2r\}}\ell(x)^3\in L^2.
\]
Hence
\[
 \mathrm{Dom}(W_rP_{\mathrm c})=\mathrm{Dom}(W_r),
\]
because
\[
 W_rP_{\mathrm c}f
 =W_rf-\langle h_\alpha,f\rangle W_rh_\alpha.
\]
Set \(B_r:=W_rP_{\mathrm c}\) on \(\mathrm{Dom}(W_r)\).  The identity
\[
 B_rf=W_rf-\langle h_\alpha,f\rangle W_rh_\alpha,
 \qquad f\in\mathrm{Dom}(W_r),
\]
shows that \(B_r\) is a bounded rank-one perturbation of a closed, densely
defined multiplication operator.  Therefore \(B_r\) is closed and densely
defined, and
\[
 B_r^*=P_{\mathrm c}W_r,\qquad
 \mathrm{Dom}(B_r^*)=\mathrm{Dom}(W_r).
\]
We next verify that this bounded form extension has the required
closed-operator meaning.  For \(\Psi\in\mathrm{Dom}(B_r^*)\), set
\[
 h=(|D|_\alpha-z)^{-1}B_r^*\Psi.
\]
For every \(f\in\mathrm{Dom}(B_r^*)\), the definition of \(T_r(z)\) gives
\[
 \langle B_r^*f,h\rangle=\langle f,T_r(z)\Psi\rangle.
\]
By the definition of the adjoint, \(h\in\mathrm{Dom}(B_r^{**})\) and
\(B_r^{**}h=T_r(z)\Psi\).  Since \(B_r\) is closed, this proves
\[
 (|D|_\alpha-z)^{-1}B_r^*\Psi\in\mathrm{Dom}(B_r)
\]
and
\[
 B_r(|D|_\alpha-z)^{-1}B_r^*\Psi=T_r(z)\Psi.
\]
In particular, \eqref{eq:cubic-epsilon}, applied at \(\bar z\), gives
\[
 \left|\left\langle (|D|_\alpha-z)^{-1}B_r^*\Psi,
 B_r^*\Psi\right\rangle\right|
 =\left|\left\langle B_r^*\Psi,
 (|D|_\alpha-\bar z)^{-1}B_r^*\Psi\right\rangle\right|
 \leq\varepsilon_r\norm \Psi_2^2.
\]
By \cite[Definition~2.1(ii) and Remark~2.1]{DAnconaKato}, this is precisely
the uniform resolvent criterion for \(|D|_\alpha\)-supersmoothness of
\(B_r\).  Kato's theory \cite{Kato}, in the form of
\cite[Theorems~2.2 and~2.3]{DAnconaKato}, gives
\eqref{eq:cubic-hom-smooth} and \eqref{eq:cubic-ret-smooth}.
\end{proof}

\subsection{The weighted Osgood argument}

We split the cubic difference into singular and regular parts.  The singular
part is controlled by localized Kato smoothing, while the regular part produces
a single logarithmic loss.  The weighted Osgood argument then yields uniqueness.
Set \(F(u):=|u|^2u\), and first prove the Duhamel formula for weak energy solutions.

\begin{lemma}\label{lem:cubic-Duhamel}
Let \(u\) be a weak energy solution on an interval \(J\) containing \(0\), and
set \(u_0:=u(0)\).  Then, for every \(t\in J\), the following identity holds in
\(L^2(\R)\)
\begin{equation}\label{eq:cubic-Duhamel}
 u(t)=e^{-it|D|_\alpha}u_0
 +i\int_0^t e^{-i(t-\tau)|D|_\alpha}F(u(\tau))\,\dd\tau.
\end{equation}
\end{lemma}

\begin{proof}
Lemma~\ref{lem:critical-form-coordinates} gives
\[
 \norm{F(u)}_2=\norm u_6^3\lesssim\norm u_{X_\alpha}^3.
\]
Strong \(L^2\)-continuity, interpolation with a uniform higher
finite exponent bound and the cubic difference estimate show that
\(F(u)\in C_{\mathrm{loc}}(L^2)\).  Let
\(P_N=\mathbf 1_{[-N,N]}(|D|_\alpha)\).  Testing the weak equation against
\(P_N\)-localized vectors gives the variation-of-constants formula for
\(P_Nu\) in \(L^2\).  Since \(P_N\to I\) strongly on \(L^2\), unitarity and
dominated convergence pass both terms to the limit uniformly on compact time
intervals, proving \eqref{eq:cubic-Duhamel}.
\end{proof}

Let \(u,v\) be weak energy solutions on \(J\) with
\(u(t_0)=v(t_0)\) for some \(t_0\in J\).  After a time translation and
restriction to a compact interval, it suffices to work on \([0,T]\), \(T>0\),
with \(u(0)=v(0)\).  Choose \(R>0\) such that
\begin{equation}\label{eq:cubic-form-ball}
 \sup_{t\in[0,T]}\bigl(\norm{u(t)}_{X_\alpha}
 +\norm{v(t)}_{X_\alpha}\bigr)\leq R.
\end{equation}
Put \(w=u-v\), \(y(t)=\norm{w(t)}_2^2\) and use the decompositions
\[
 u=\phi_u+q_uG_{2\lambda_\alpha},\qquad
 v=\phi_v+q_vG_{2\lambda_\alpha}.
\]
Then
\begin{equation}\label{eq:cubic-decomposition-bound}
 \sup_{t\in[0,T]}\bigl(\norm{\phi_u(t)}_{H^{1/2}}
 +\norm{\phi_v(t)}_{H^{1/2}}+|q_u(t)|+|q_v(t)|\bigr)
 \leq C_\alpha R.
\end{equation}
Choose \(K>1\), depending on \(R\) and \(\alpha\), such that for every \(t\in[0,T]\)
\[
 0\leq y(t)\leq K/e,\qquad \norm{w(t)}_4^2\leq K.
\]
Set \(L_K(\rho):=1+\log(K/\rho)\) for \(\rho>0\), and extend
\(\rho L_K(\rho)\) and \(\sqrt{\rho}L_K(\rho)\) continuously by zero at
\(\rho=0\).
We now decompose the cubic difference into singular and regular parts.

\begin{lemma}
\label{lem:cubic-products}
Let \(u,v\) be two weak energy solutions on \([0,T]\) with the same initial
datum. Then for every fixed \(0<r<r_0\), there exist measurable functions
\(F_{\mathrm{sing}}\) and \(F_{\mathrm{reg}}:\R\times\R\to\C\), such that the following relations hold
for every \(t\in[0,T]\):
\begin{enumerate}[\upshape(i)]
\item
\[
 F(u)-F(v)=W_rF_{\mathrm{sing}}+F_{\mathrm{reg}}.
\]
\item
\[
 \norm{F_{\mathrm{sing}}(t)}_2
 \leq C_{R,r}\sqrt{y(t)}.
\]
\item
\[
 \norm{F_{\mathrm{reg}}(t)}_2
 \leq C_{R,r}\sqrt{y(t)}L_K(y(t)).
\]
\end{enumerate}
\end{lemma}

\begin{proof}
The real derivative of \(F\) gives the exact representation
\[
 F(u)-F(v)=a_1(u,v)w+a_2(u,v)\overline w,
\]
where, with \(z_s=v+sw\),
\[
 a_1=2\int_0^1|z_s|^2\,\dd s,
 \qquad
 a_2=\int_0^1z_s^2\,\dd s.
\]
By \eqref{eq:cubic-decomposition-bound} and the elementary cubic inequality, the
Green contribution on \(\{|x|<r\}\) is bounded by a multiple of
\(W_r\).  On the cutoff annulus and on \(\{|x|\geq r\}\), the Green
function is bounded, so that contribution is absorbed into the constant term
below.  Hence
\[
 |a_1|+|a_2|\leq S+R_{\mathrm{reg}},
\]
where
\[
 S=C_RW_r,\qquad
 R_{\mathrm{reg}}=C_{R,r}\bigl(1+|\phi_u|^2+|\phi_v|^2\bigr).
\]
For \(j=1,2\), define
\[
 a_{j,\mathrm{sing}}=a_j\frac{S}{S+R_{\mathrm{reg}}},\qquad
 a_{j,\mathrm{reg}}=a_j\frac{R_{\mathrm{reg}}}{S+R_{\mathrm{reg}}},
\]
with both quotients set to zero where \(S+R_{\mathrm{reg}}=0\).  This gives an exact
measurable splitting with \(|a_{j,\mathrm{sing}}|\leq S\) and
\(|a_{j,\mathrm{reg}}|\leq R_{\mathrm{reg}}\).  Set
\[
 W_rF_{\mathrm{sing}}
 =a_{1,\mathrm{sing}}w+a_{2,\mathrm{sing}}\overline w
\]
and
\[
 F_{\mathrm{reg}}=a_{1,\mathrm{reg}}w+a_{2,\mathrm{reg}}\overline w,
\]
where \(F_{\mathrm{sing}}=0\) on \(\{W_r=0\}\).  Then
\(|F_{\mathrm{sing}}|\leq C_R|w|\), which gives the first estimate.

It remains to estimate \(|\phi|^2w\).  For \(p\geq2\), H\"older's inequality and
\eqref{eq:critical-regular-Lp} give
\[
 \norm{|\phi|^2w}_2
 \leq C_Rp
       \norm w_{2p/(p-1)}.
\]
Interpolation between \(L^2\) and the uniform \(L^4\)-bound gives
\[
 \norm w_{2p/(p-1)}
 \leq C_R y^{1/2-1/p}K^{1/p}.
\]
If \(y=0\), the estimate is immediate.  For \(y>0\), choose
\(p=L_K(y)\), which is at least \(2\).  Since
\((K/y)^{1/p}\leq e\), the second estimate follows.
\end{proof}

Define
\[
 \mathcal Q(t)=\int_0^t y(\tau)\,\dd\tau,
\]
\[
 \mathcal Y(t)=\int_0^t y(\tau)L_K(y(\tau))\,\dd\tau,
\]
\[
 \mathcal J(t)=\int_0^t\sqrt{y(\tau)}
 L_K(y(\tau))\,\dd\tau.
\]
\begin{proposition}
\label{prop:cubic-weighted-difference}
For every fixed \(0<r<r_0\),
\begin{equation}\label{eq:cubic-y-bound}
 y(t)\leq C_{R,r,T}
 \bigl(\mathcal Y(t)+\mathcal J(t)^2\bigr),
 \qquad t\in[0,T].
\end{equation}
\end{proposition}

\begin{proof}
We first establish an auxiliary weighted spacetime estimate.  Set
\[
 \mathcal X(t):=\norm{W_rw}_{L^2([0,t]\times\R)}.
\]
This quantity is finite since
\[
 \norm{W_r\phi}_2\leq\norm{W_r}_4\norm\phi_4,
 \qquad \norm{W_rG_{2\lambda_\alpha}}_2<\infty.
\]
We first control the discrete component.  Let
\(b(t)=\langle h_\alpha,w(t)\rangle\).  Projecting the difference equation
onto \(h_\alpha\), the singular source is bounded by
\(\norm{W_rh_\alpha}_2\norm{F_{\mathrm{sing}}}_2\), and the regular source by
\(\norm{F_{\mathrm{reg}}}_2\).  Since \(L_K\geq1\) and \(b(0)=0\),
Lemma~\ref{lem:cubic-products} gives
\begin{equation}\label{eq:cubic-bound-state}
 |b(t)|\leq C_{R,r}\mathcal J(t).
\end{equation}
Since \(W_rh_\alpha\in L^2\) and
\((I-P_{\mathrm c})w=b(t)h_\alpha\), it follows that
\begin{equation}\label{eq:cubic-bound-state-weighted}
 \norm{W_r(I-P_{\mathrm c})w}_{L^2([0,t]\times\R)}
 \leq C_{R,r,T}\mathcal J(t).
\end{equation}

For the continuous component, work on \(L^2(\R)\oplus L^2(\R)\) with
\[
 \mathbf w=(w,\overline w)^{\mathsf T},\qquad
 \mathbf A=\operatorname{diag}(|D|_\alpha,-|D|_\alpha),\qquad
 \mathbf P_{\mathrm c}=\operatorname{diag}(P_{\mathrm c},P_{\mathrm c}).
\]
Let \(\mathcal F=F(u)-F(v)\) and
\(\mathbf F=(\mathcal F,-\overline{\mathcal F})^{\mathsf T}\).  Since
\(\mathbf w(0)=0\),
Lemma~\ref{lem:cubic-Duhamel} gives
\begin{equation}\label{eq:cubic-projected-doubled-Duhamel}
 \mathbf P_{\mathrm c}\mathbf w(\tau)
 =i\int_0^\tau e^{-i(\tau-\sigma)\mathbf A}
 \mathbf P_{\mathrm c}\mathbf F(\sigma)\,\dd\sigma.
\end{equation}
For the singular source in Lemma~\ref{lem:cubic-products},
\[
 \norm{F_{\mathrm{sing}}}_{L^2([0,t]\times\R)}
 \leq C_{R,r}\mathcal Q(t)^{1/2}.
\]
We now identify the Kato extension with the corresponding Bochner Duhamel
term.  The pointwise cubic bound
\[
 |W_rF_{\mathrm{sing}}|
 \leq C_R\bigl(|u|^2+|v|^2\bigr)|w|
\]
together with the finite exponent embeddings implies
\[
 F_{\mathrm{sing}}\in L^2([0,t];\mathrm{Dom}(W_r)).
\]
Indeed, Lemma~\ref{lem:cubic-products} controls its \(L^2\)-part, while
H\"older's inequality gives
\[
 \norm{W_rF_{\mathrm{sing}}(\tau)}_2
 \leq C_R\bigl(\norm{u(\tau)}_8^2+\norm{v(\tau)}_8^2\bigr)
 \norm{w(\tau)}_4.
\]
The right-hand side is uniformly bounded on \([0,t]\).  The graph space
\(\mathrm{Dom}(W_r)\) is separable, so its simple functions are dense in the
corresponding Bochner space.
Extend \(F_{\mathrm{sing}}\) by zero outside \([0,t]\).  There are
\(\mathrm{Dom}(W_r)\)-valued step functions \(F_m\), supported in \([0,t]\),
such that
\[
 F_m\longrightarrow F_{\mathrm{sing}}
 \quad\hbox{in }L^2(\R;\mathrm{Dom}(W_r))
\]
for the graph norm of \(W_r\).  For \(0\leq\tau\leq t\), set
\[
 U_m(\tau)=\int_0^\tau e^{-i(\tau-\sigma)|D|_\alpha}
 P_{\mathrm c}W_rF_m(\sigma)\,\dd\sigma
\]
and
\[
 U(\tau)=\int_0^\tau e^{-i(\tau-\sigma)|D|_\alpha}
 P_{\mathrm c}W_rF_{\mathrm{sing}}(\sigma)\,\dd\sigma.
\]
Both are ordinary Bochner integrals in \(L^2\).  For every \(m\),
\cite[Theorem~2.3]{DAnconaKato} identifies the Kato extension with
\(B_rU_m=W_rU_m\), since \(P_{\mathrm c}U_m=U_m\).  Proposition~\ref{prop:cubic-Kato}
shows that \(B_rU_m\) converges in \(L^2([0,t]\times\R)\) to some function
\(V\).  On the other hand, unitarity gives
\[
 \sup_{0\leq\tau\leq t}\norm{U_m(\tau)-U(\tau)}_2
 \leq t^{1/2}\norm{W_r(F_m-F_{\mathrm{sing}})}_{L^2([0,t]\times\R)},
\]
which tends to zero.  After passing to a subsequence,
\(B_rU_m(\tau)\to V(\tau)\) in \(L^2\) for almost every \(\tau\).  The
closedness of \(B_r\) then gives \(U(\tau)\in\mathrm{Dom}(B_r)\) and
\(B_rU(\tau)=V(\tau)\) for almost every \(\tau\).  Thus the Kato extension is
exactly \(B_rU=W_rU\), where \(U\) is the ordinary Bochner Duhamel term for
the singular source in \eqref{eq:cubic-projected-doubled-Duhamel}.

The zero extension and the restriction back to \([0,t]\) do not increase the
spacetime norm.  The same argument applies to the conjugate component because
\((-|D|_\alpha-z)^{-1}=-(|D|_\alpha+z)^{-1}\), and the direct-sum estimate
follows componentwise with the same constant.  The retarded estimate in
Proposition~\ref{prop:cubic-Kato} therefore gives
\[
 \left\|W_r\int_0^\tau e^{-i(\tau-\sigma)|D|_\alpha}
 P_{\mathrm c}W_rF_{\mathrm{sing}}(\sigma)\,\dd\sigma
 \right\|_{L^2([0,t]_\tau\times\R)}
 \leq C_{\alpha,r}\norm{F_{\mathrm{sing}}}_{L^2([0,t]\times\R)}
 \leq C_{R,r}\mathcal Q(t)^{1/2}.
\]

For the regular source, the homogeneous estimate followed by Minkowski's
inequality gives
\[
 \left\|W_r\int_0^\tau e^{-i(\tau-\sigma)|D|_\alpha}
 P_{\mathrm c}F_{\mathrm{reg}}(\sigma)\,\dd\sigma
 \right\|_{L^2([0,t]_\tau\times\R)}
 \leq C_{R,r}\mathcal J(t).
\]
Combining the continuous estimates with the discrete estimate
\eqref{eq:cubic-bound-state-weighted}, and using
\(\mathcal Q\leq\mathcal Y\), we obtain
\begin{equation}\label{eq:cubic-X-bound}
\mathcal X(t)\leq C_{R,r,T}
 \bigl(\mathcal Y(t)^{1/2}+\mathcal J(t)\bigr).
\end{equation}

We now derive the stated \(L^2\)-distance estimate.
Time mollification and the Hilbert triple chain rule give
\(y\in W^{1,\infty}([0,T])\) and
\(y'=2\operatorname{Re}\langle\partial_tw,w\rangle_{X_\alpha^*,X_\alpha}\)
almost everywhere.  Indeed, the identity holds for the time-mollified
functions, and convergence in the Hilbert triple passes it to the limit.
Self-adjointness cancels the linear form.  Pairing the singular part with
\(w\), the regular part with \(w\) in \(L^2\) and integrating in time gives
\[
 y(t)\leq C_{R,r}
 \bigl(\mathcal X(t)\mathcal Q(t)^{1/2}+\mathcal Y(t)\bigr).
\]
Insert \eqref{eq:cubic-X-bound}, use
\(\mathcal Q\leq\mathcal Y\) and apply Young's inequality to
\(\mathcal J\mathcal Y^{1/2}\).  This proves
\eqref{eq:cubic-y-bound}.
\end{proof}

The preceding distance estimate reduces uniqueness to an Osgood argument.
We now complete the proof of Proposition~\ref{prop:cubic-half-wave-weak-unique}.

\begin{proof}[Proof of Proposition~\ref{prop:cubic-half-wave-weak-unique}]
On \([0,T]\), \eqref{eq:cubic-form-ball} holds and
Proposition~\ref{prop:cubic-weighted-difference} applies.
By definition, \(\mathcal Y\) and \(\mathcal J\) are nonnegative and
absolutely continuous, vanish at zero and satisfy
\[
 \mathcal Y'=yL_K(y)
\]
and
\[
 \mathcal J'=\sqrt y\,L_K(y)
\]
almost everywhere.

Set \(\Omega=\mathcal Y+\mathcal J^2\), and work up to the first time at which \(\Omega\)
reaches a sufficiently small fixed number.  On this interval the function
\[
 f(\rho)=\rho L_K(\rho)
\]
is increasing.  For \(a,b\) in the same range,
\begin{equation}\label{eq:cubic-log-product}
 \sqrt{ab}\,L_K(b)
 \leq C\bigl(f(a)+f(b)\bigr).
\end{equation}
Indeed, if \(b\geq a\), the left side is bounded by \(f(b)\).  If
\(b<a\), write \(b=\theta a\); then
\(L_K(b)=L_K(a)+\log(1/\theta)\) and use
\(\sup_{0<\theta<1}\sqrt\theta\,\log(1/\theta)<\infty\).

Apply \eqref{eq:cubic-log-product} with \(a=\mathcal J^2\), \(b=y\).  Since
\(y,\mathcal J^2\leq C\Omega\), monotonicity gives
\[
 \Omega'\leq C f(C\Omega)
 \leq C\Omega L_K(\Omega).
\]
The Osgood integral
\[
 \int_{0^+}\frac{\dd\rho}{\rho L_K(\rho)}
\]
diverges.  Osgood's lemma yields
\(\Omega=0\) before the first exit, so the exit cannot occur.  Hence
\eqref{eq:cubic-y-bound} gives \(y=0\) on \([0,T]\).
Time reversal and complex conjugation give the backward conclusion, and
compact intervals exhaust the common lifespan.
\end{proof}
The preceding proof shows that the present Osgood argument for the mass-critical (cubic) problem can be extended to mass-subcritical ones. For mass-supercritical cases, a stronger logarithmic loss would make the corresponding Osgood integral finite, so the approach does not work. In this sense, the mass-critical exponent is exactly the threshold for the present Osgood argument, which coincides with the corresponding result proved by Krieger, Lenzmann and Rapha\"{e}l \cite{KriegerLenzmannRaphael} for the classical free half-wave equation. Although there are some well-posedness results for mass-supercritical problems in some spaces with higher regularity, such as the results od Ozawa and Visciglia \cite{OzawaVisciglia}, unconditional uniqueness for mass-supercritical problems in energy space is still open.

\section{Well-posedness for the Cauchy problem}
\label{sec:cubic-flow}

In this section, we first combine an abstract energy method with the unconditional uniqueness
established in Section~\ref{sec:cubic-linear} to obtain local well-posedness in
\(X_\alpha\) and the conservation laws.  We then establish a coercive bound
below the critical mass.  Together with the uniform local estimates, these
ingredients prove Theorem~\ref{thm:intro-cubic-flow}.

\subsection{Local well-posedness}

We will apply an abstract theory of existence and uniqueness given by Okazawa, Suzuki and Yokota
\cite{OkazawaSuzukiYokota}.

Let \(S\) be a nonnegative selfadjoint operator in a complex Hilbert space
\(X\).  Put \(X_S:=D((1+S)^{1/2})\).  Then we have the usual triplet:
\(X_S\subset X=X^*\subset X_S^*\).  Under this setting \(S\) can be extended
to a nonnegative selfadjoint operator in \(X_S^*\) with domain \(X_S\). Consider the following abstract Cauchy problem:
\begin{equation*}
 \left\{
 \begin{array}{l}
  i\dfrac{\dd u}{\dd t}=Su+g(u),\\[4pt]
  u(0)=u_0,
 \end{array}
 \right.
 \tag{ACP}
\end{equation*}
where \(g:X_S\to X_S^*\) is a nonlinear operator.  

\medskip
\noindent\textbf{(G1) Existence of energy functional:}
there exists \(G\in C^1(X_S;\R)\) such that \(G'=g\), that is, given
\(u\in X_S\), for every \(\varepsilon>0\) there exists
\(\delta=\delta(u,\varepsilon)>0\) such that
\[
 \abs{G(u+v)-G(u)-\operatorname{Re}\langle g(u),v\rangle_{X_S^*,X_S}}
 \leq\varepsilon\norm v_{X_S}
\]
for all \(v\in X_S\) with \(\norm v_{X_S}<\delta\).

\medskip
\noindent\textbf{(G2) Local Lipschitz continuity:}
for all \(M>0\) there exists \(C(M)>0\) such that
\[
 \norm{g(u)-g(v)}_{X_S^*}
 \leq C(M)\norm{u-v}_{X_S}
\]
for all \(u,v\in X_S\) with
\(\norm u_{X_S},\norm v_{X_S}\leq M\).

\medskip
\noindent\textbf{(G3) H\"older-like continuity of energy functional:}
given \(M>0\), for all \(\delta>0\) there exists a constant
\(C_\delta(M)>0\) such that
\[
 \abs{G(u)-G(v)}
 \leq\delta+C_\delta(M)\norm{u-v}_X
\]
for all \(u,v\in X_S\) with
\(\norm u_{X_S},\norm v_{X_S}\leq M\).

\medskip
\noindent\textbf{(G4) Gauge type condition:}
\[
 \operatorname{Re}\langle g(u),iu\rangle_{X_S^*,X_S}=0
 \qquad\text{for all }u\in X_S.
\]

\medskip
\noindent\textbf{(G5) Closedness type condition:}
given a bounded open interval \(I\subset\R\), let \((w_n)\) be any bounded
sequence in \(L^\infty(I;X_S)\) such that
\begin{equation*}
 \left\{
 \begin{array}{ll}
  w_n(t)\rightharpoonup w(t) & \text{weakly in }X_S
  \text{ for almost every }t\in I,\\
  g(w_n)\overset{*}{\rightharpoonup}f
  & \text{in }L^\infty(I;X_S^*).
 \end{array}
 \right.
\end{equation*}
Then
\begin{equation*}
 \int_I \operatorname{Re}\langle f(t),iw(t)\rangle_{X_S^*,X_S}\,\dd t
 =\lim_{n\to\infty}\int_I
  \operatorname{Re}\langle g(w_n(t)),iw_n(t)\rangle_{X_S^*,X_S}\,\dd t.
\end{equation*}
Here \(f=g(w)\) is guaranteed if
\begin{equation*}
 w_n(t)\to w(t)\quad\text{strongly in }X
 \quad\text{for almost every }t\in I.
\end{equation*}

Okazawa, Suzuki and Yokota proved that, under assumptions
\(\mathrm{(G1)}\)--\(\mathrm{(G5)}\), \(\mathrm{(ACP)}\) admits a local weak
solution with conserved \(X\)-norm and non-increasing energy \cite[Theorem~2.2]{OkazawaSuzukiYokota}. They also point out that if the uniqueness of local weak solutions is further assumed, then
the solution enjoys stronger temporal regularity, and both its \(X\)-norm and
energy are conserved \cite[Theorem~2.3]{OkazawaSuzukiYokota}.  Applying these
results in the present setting yields the following proposition:

\begin{proposition}
\label{prop:cubic-half-wave-local-flow}
For every \(u_0\in X_\alpha\), there is \(T>0\) such that
\eqref{eq:evolution} has a unique solution
\[
 u\in C((-T,T);X_\alpha)\cap C^1((-T,T);X_\alpha^*).
\]
Moreover, for every \(t\in(-T,T)\),
\[
 M(u(t))=M(u_0),
 \qquad
 E_\alpha(u(t))=E_\alpha(u_0).
\]
\end{proposition}

\begin{proof}
We apply the abstract results in \cite{OkazawaSuzukiYokota} with
\[
 S=|D|_\alpha+\lambda_\alpha.
\]
\[
 g(u)=-\lambda_\alpha u-F(u).
\]
Here \(S\) is a nonnegative self-adjoint operator in \(L^2(\R)\), with
\(\operatorname{Dom}(S^{1/2})=X_\alpha\).  Moreover,
\(\norm u_{X_\alpha}^2=\norm{S^{1/2}u}_2^2+\lambda_\alpha\norm u_2^2\),
so this norm is equivalent to the usual graph norm.  Equation
\eqref{eq:evolution} becomes
\[
 i\partial_tu=Su+g(u).
\]
The associated real potential is
\[
 P_\alpha(u)
 =-\frac{\lambda_\alpha}{2}\norm u_2^2
  -\frac14\norm u_4^4,
\]
and
\[
 E_\alpha(u)
 =\frac12\norm{S^{1/2}u}_2^2+P_\alpha(u).
\]
The derivative of \(P_\alpha\) is
\[
 DP_\alpha(u)[h]
 =\operatorname{Re}\langle g(u),h\rangle_{X_\alpha^*,X_\alpha}.
\]
We first verify \(\mathrm{(G1)}\)--\(\mathrm{(G5)}\).

\smallskip
\noindent\emph{Verification of \(\mathrm{(G1)}\).}
Lemma~\ref{lem:critical-form-coordinates} shows that \(P_\alpha\) is
well-defined on \(X_\alpha\).  The derivative formula above and the cubic
difference estimate below show that its derivative \(g\) is continuous from
\(X_\alpha\) to \(X_\alpha^*\).  Thus \(P_\alpha\in C^1(X_\alpha;\R)\), which
proves \(\mathrm{(G1)}\).

\smallskip
\noindent\emph{Verification of \(\mathrm{(G2)}\).}
The pointwise estimate
\[
 \norm{F(u)-F(v)}_2
 \leq C\bigl(\norm u_6^2+\norm v_6^2\bigr)\norm{u-v}_6,
\]
together with \(X_\alpha\hookrightarrow L^6\) and
\(L^2\hookrightarrow X_\alpha^*\), gives \(\mathrm{(G2)}\).

\smallskip
\noindent\emph{Verification of \(\mathrm{(G3)}\).}
On bounded subsets of \(X_\alpha\),
\[
 \abs{\norm u_4^4-\norm v_4^4}
 \leq C\bigl(\norm u_6^3+\norm v_6^3\bigr)\norm{u-v}_2,
\]
and the analogous quadratic estimate for the mass gives
\[
 \abs{P_\alpha(u)-P_\alpha(v)}
 \leq C_R\norm{u-v}_2.
\]
This proves the stronger form of \(\mathrm{(G3)}\) with \(\delta=0\).

\smallskip
\noindent\emph{Verification of \(\mathrm{(G4)}\).}
Gauge covariance gives \(\mathrm{(G4)}\); under our pairing convention this is
equivalent to \(\operatorname{Im}\langle g(u),u\rangle=0\).

\smallskip
\noindent\emph{Verification of \(\mathrm{(G5)}\).}
Under its hypotheses, the local compactness
in Lemma~\ref{lem:critical-form-coordinates} and uniqueness of the weak limit
give strong local \(L^6\)-convergence at almost every time.  Boundedness in
\(L^\infty(J;X_\alpha)\) and the embedding
\(X_\alpha\hookrightarrow L^6(\R)\) give a uniform bound in
\(L^\infty(J;L^6(\R))\).  Dominated convergence then implies that
\(F(u_n)\to F(u)\) strongly on compact spacetime cylinders.  Moreover,
\[
 \norm{F(u_n)}_{L^2(J;L^2)}
 \leq |J|^{1/2}\norm{u_n}_{L^\infty(J;L^6)}^3,
\]
so \((F(u_n))\) is bounded in \(L^2(J;L^2)\).  Testing against compactly
supported functions and then using their density identifies every global weak
cluster point and yields
\begin{equation}\label{eq:critical-nonlinear-spacetime-weak}
 F(u_n)\rightharpoonup F(u)
 \quad\text{in }L^2(J;L^2(\R)).
\end{equation}
The embedding \(L^2\hookrightarrow X_\alpha^*\) gives the same closure in
\(L^2(J;X_\alpha^*)\), also for
\(g(v)=-\lambda_\alpha v-F(v)\).  Hence the weak-star limit in
\(\mathrm{(G5)}\) is \(f=g(u)\).  By \(\mathrm{(G4)}\), both sides of the
identity in \(\mathrm{(G5)}\) vanish.  Thus \(\mathrm{(G1)}\)--\(\mathrm{(G5)}\)
hold.  Theorem~2.2 in \cite{OkazawaSuzukiYokota} therefore gives a weak energy
solution with conserved mass and the energy inequality.
Proposition~\ref{prop:cubic-half-wave-weak-unique} gives uniqueness in the weak
energy class.  Hence \cite[Theorem~2.3]{OkazawaSuzukiYokota} yields the stated
time regularity and the conservation of mass and energy.
\end{proof}

\subsection{Global continuation below the critical mass}

The local flow and Lemma~\ref{lem:cubic-half-wave-coercivity} now give the
maximal dynamics.

\begin{proof}[Proof of Theorem~\ref{thm:intro-cubic-flow}]
Proposition~\ref{prop:cubic-half-wave-local-flow} and uniqueness on
overlapping intervals produce a unique maximal solution on
\((-T_-,T_+)\).  The conservation laws hold throughout this interval.

The construction is uniform on bounded subsets of \(X_\alpha\).  More precisely,
\cite[Theorem~2.2, Lemma~4.2 and Eqs.~(4.12)--(4.13)]{OkazawaSuzukiYokota}
gives, for every \(R>0\), constants \(T_R,C_R>0\) such that all data with
\(\norm{v_0}_{X_\alpha}\leq R\) admit solutions on \((-T_R,T_R)\) satisfying
\[
 \norm{v}_{L^\infty((-T_R,T_R);X_\alpha)}
 +\norm{\partial_tv}_{L^\infty((-T_R,T_R);X_\alpha^*)}
 \leq C_R.
\]
The bounds for the approximating solutions pass to the weak limit by lower
semicontinuity.

We next prove continuous dependence.  Let \(u_{0,n}\to u_0\) in
\(X_\alpha\), and let \(u_n\) and \(u\) be the corresponding maximal
solutions.  The preceding uniform estimate gives a common compact local
interval \(J\) containing zero and
\[
 \sup_n\left(
 \norm{u_n}_{L^\infty(J;X_\alpha)}
 +\norm{\partial_tu_n}_{L^\infty(J;X_\alpha^*)}
 \right)<\infty.
\]
The weak Arzel\`a--Ascoli theorem
\cite[Proposition~1.1.2(i)]{Cazenave} gives, after extraction from any
subsequence, \(u_n\to\widetilde u\) in \(C_w(J;X_\alpha)\).  The local
compactness and cubic convergence argument in the proof of
Proposition~\ref{prop:cubic-half-wave-local-flow} identifies the nonlinear
limit.  Passing to the weak equation shows that \(\widetilde u\) is a weak
energy solution with initial datum \(u_0\).  Unconditional uniqueness gives
\(\widetilde u=u\).  Thus the full sequence converges in
\(C_w(J;X_\alpha)\).

The interpolation inequality
\(\norm f_2^2\lesssim\norm f_{X_\alpha}\norm f_{X_\alpha^*}\) gives
equicontinuity in \(L^2\).  If \(t_n\to t\), every weak \(L^2\) cluster point
of \(u_n(t_n)\) is \(u(t)\).  Since
\(M(u_n(t_n))=M(u_{0,n})\to M(u_0)=M(u(t))\), weak convergence and convergence
of the norms imply strong convergence.  Compactness of \(J\) then gives
\[
 u_n\to u\qquad\text{in }C(J;L^2).
\]
Interpolation yields convergence in \(C(J;L^4)\).  For every solution \(v\)
with initial datum \(v_0\), the conservation laws give
\begin{equation}\label{eq:cubic-form-norm-flow-identity}
 \norm{v(t)}_{X_\alpha}^2
 =2E_\alpha(v_0)+\frac12\norm{v(t)}_4^4
  +2\lambda_\alpha\norm{v_0}_2^2.
\end{equation}
For any \(t_n\to t\), the uniform \(X_\alpha\)-bound and the preceding
\(L^2\)-convergence identify every weak \(X_\alpha\) cluster point of
\(u_n(t_n)\) with \(u(t)\).  Equation
\eqref{eq:cubic-form-norm-flow-identity} and
\(E_\alpha(u_{0,n})\to E_\alpha(u_0)\) give
\(\norm{u_n(t_n)}_{X_\alpha}\to\norm{u(t)}_{X_\alpha}\).  Uniform convexity
and the same sequential argument therefore yield
\[
 u_n\to u\qquad\text{in }C(J;X_\alpha).
\]

Let \(K\subset(-T_-,T_+)\) be compact.  Enlarge it, if necessary, to a compact
interval containing zero.  The image of this interval under \(u\) is bounded
in \(X_\alpha\).  The uniform local existence time on bounded sets and the
local convergence just proved allow a finite sequence of restarts along the
interval.  Hence \(K\) lies in the lifespan of \(u_n\) for all large \(n\),
and \(u_n\to u\) in \(C(K;X_\alpha)\).  This proves the continuous dependence
asserted in part~\(\mathrm{(i)}\).

Suppose that \(T_+<\infty\) and the asserted limit fails.  Then there are
\(R>0\) and \(t_n\uparrow T_+\) such that
\(\norm{u(t_n)}_{X_\alpha}\leq R\).  The uniform local existence time above
allows us to restart the
solution at each \(t_n\).  For all sufficiently large \(n\),
\(t_n+T_R>T_+\).  Uniqueness identifies the restarted solution with \(u\) on
their common interval.  This extends \(u\) beyond \(T_+\), which contradicts
maximality.  The argument at \(-T_-\) is identical.  Therefore, whenever
\(T_\pm<\infty\),
\[
 \lim_{t\to\pm T_\pm}\norm{u(t)}_{X_\alpha}=\infty.
\]

Finally, if \(0<M(u_0)<\mu_*\), Lemma~\ref{lem:cubic-half-wave-coercivity}
and the conservation laws bound the solution in \(X_\alpha\) throughout its
maximal lifespan.  Consequently,
\[
 T_\pm=\infty,
 \qquad
 \sup_{t\in\R}\norm{u(t)}_{X_\alpha}<\infty.
\]
The zero solution covers \(M(u_0)=0\).
\end{proof}

\section{Ground states and orbital stability}
\label{sec:cubic-ground-states}

In this section we first show that the constrained energy is unbounded from
below at and above the critical mass, and then prove the existence and orbital
stability of ground states below the critical mass.

\subsection{Constrained energy levels}

The analysis of the constrained energy at the critical mass is based on the
following Green perturbation of the free optimizer. Recall that \(Q_*\) is the unique positive and even ground state of \eqref{eq:intro-free-ground-state}.

\begin{lemma}
\label{lem:cubic-critical-green-variation}
Fix \(\lambda>\lambda_\alpha\).  For sufficiently small \(q>0\), mass
normalization of the Green perturbation \(Q_*+qG_\lambda\) gives a family
\(v_q\in X_\alpha\) with \(M(v_q)=M(Q_*)=\mu_*\).  As \(q\to0\), there is
\(\phi_q\in H^{1/2}(\R)\) with
\(\norm{\phi_q-Q_*}_{H^{1/2}}=O(q)\) such that
\begin{equation}\label{eq:cubic-critical-charge-derivative}
 v_q=\phi_q+\bigl(q+O(q^2)\bigr)G_\lambda,
 \qquad
 E_\alpha(v_q)=-Q_*(0)q+o(q).
\end{equation}
\end{lemma}

\begin{proof}
Explicitly, set
\[
 v_q:=\frac{\sqrt{\mu_*}(Q_*+qG_\lambda)}
 {\norm{Q_*+qG_\lambda}_2}.
\]
Then \(M(v_q)=M(Q_*)=\mu_*\).  Extend this curve to real \(q\) near zero and write
\(v_q=a(q)(Q_*+qG_\lambda)\), where
\(a(q)=\sqrt{\mu_*}/\norm{Q_*+qG_\lambda}_2\), and let
\(c:=\langle Q_*,G_\lambda\rangle\in\R\).  One has
\[
 a(0)=1,
 \qquad
 a'(0)=-\frac{c}{\mu_*},
\]
so \(\phi_q:=a(q)Q_*\) satisfies
\(\norm{\phi_q-Q_*}_{H^{1/2}}=O(q)\), while the Green coefficient of
\(v_q\) is \(a(q)q=q+O(q^2)\).  Moreover,
\eqref{eq:critical-form} gives
\[
 \left.\frac{\dd}{\dd q}\right|_{q=0}
 \mathfrak t_\alpha[Q_*+qG_\lambda]=-2\lambda c.
\]
Since \(Q_*^3\in L^2(\R)\), the profile equation implies
\(Q_*\in H^1(\R)\hookrightarrow C_0(\R)\).  Pairing it with
\(G_\lambda\), and using \((|D|+\lambda)G_\lambda=\delta_0\), gives
\begin{equation}\label{eq:free-profile-green-pairing}
 \int_\R Q_*^3G_\lambda\,\dd x=Q_*(0)+(1-\lambda)c.
\end{equation}
Testing the same equation against \(Q_*\) gives
\[
 \norm{|D|^{1/2}Q_*}_2^2-\norm{Q_*}_4^4=-\mu_*.
\]
The real differentiability of the fourth-power functional and the preceding
identities yield
\[
 \left.\frac{\dd}{\dd q}\right|_{q=0}E_\alpha(v_q)
 =c-\lambda c-\int_\R Q_*^3G_\lambda\,\dd x=-Q_*(0).
\]
Since \(v_0=Q_*\) has zero charge,
\(E_\alpha(v_0)=E_\infty(Q_*)=0\), and this proves
\eqref{eq:cubic-critical-charge-derivative}.
\end{proof}

We now give the complete classification of the constrained energy levels.

\begin{proposition}
\label{prop:cubic-level-classification}
For \(0<\mu<\mu_*\), the constrained energy is finite and satisfies
\begin{equation}\label{eq:cubic-half-wave-strict-binding}
 -\infty<I_\alpha(\mu)
 <-\frac{\lambda_\alpha}{2}\mu<0.
\end{equation}
If \(\mu\geq\mu_*\), then
\[
 I_\alpha(\mu)=-\infty,
 \qquad
 \mathcal M_\alpha(\mu)=\varnothing.
\]
\end{proposition}

\begin{proof}
For \(0<\mu<\mu_*\), Lemma~\ref{lem:cubic-half-wave-coercivity} gives
finiteness from below.  Testing with \(\sqrt\mu\,h_\alpha\), where
\(M(h_\alpha)=1\) and
\(|D|_\alpha h_\alpha=-\lambda_\alpha h_\alpha\), gives
\[
 I_\alpha(\mu)
 \leq-\frac{\lambda_\alpha}{2}\mu
      -\frac{\mu^2}{4}\norm{h_\alpha}_4^4
 <-\frac{\lambda_\alpha}{2}\mu,
\]
which proves \eqref{eq:cubic-half-wave-strict-binding}.

If \(\mu>\mu_*\), set \(a=(\mu/\mu_*)^{1/2}>1\).  The computation in
Lemma~\ref{lem:free-cubic-geometry} gives
\(E_\infty(aQ_*)=\frac\mu2(1-\mu/\mu_*)<0\).  Since \(aQ_*\) has zero
charge, \(E_\alpha(aQ_*)=E_\infty(aQ_*)\).  The mass-preserving dilation
then sends the energy to \(-\infty\).

It remains to consider \(\mu=\mu_*\).  Fix \(\lambda>\lambda_\alpha\) and
let \(v_q\) be the family from
Lemma~\ref{lem:cubic-critical-green-variation}.
For \(u=\phi+qG_\lambda\in X_\alpha\),
Lemma~\ref{lem:green-coordinates} and
\(\Phi_\alpha(-\rho\lambda)
=\Phi_\alpha(-\lambda)+\pi^{-1}\log\rho\) give
\begin{equation}\label{eq:cubic-energy-scaling-anomaly}
 E_\alpha(\D_\rho u)
 =\rho E_\alpha(u)+\frac{\rho\log\rho}{2\pi}|q|^2.
\end{equation}
Choose \(0<\kappa<2\pi Q_*(0)\).  Since the charge of
\(v_{\kappa/\log\rho}\) equals
\(\kappa/\log\rho+O((\log\rho)^{-2})\), combining
\eqref{eq:cubic-critical-charge-derivative} with
\eqref{eq:cubic-energy-scaling-anomaly} gives
\[
 E_\alpha(\D_\rho v_{\kappa/\log\rho})
 =\frac{\rho}{\log\rho}
 \left(-\kappa Q_*(0)+\frac{\kappa^2}{2\pi}+o(1)\right)
 \to-\infty.
\]
Thus the level is \(-\infty\) at and above \(\mu_*\), which also rules out
minimizers.
\end{proof}

\begin{remark}
The resolvent formula in Lemma~\ref{lem:operator-regimes} shows that
\(|D|_\alpha\) converges to the free half-wave operator \(|D|\) in the norm
resolvent sense as \(\alpha\to+\infty\).  The critical fixed-mass variational
problem nevertheless has a singular limit.  Proposition~\ref{prop:cubic-level-classification}
gives \(I_\alpha(\mu_*)=-\infty\) for every finite \(\alpha\), whereas
Lemma~\ref{lem:free-cubic-geometry} gives \(I_\infty(\mu_*)=0\), with the
infimum attained by the free optimizers.  Thus linear norm resolvent
convergence does not imply convergence of the critical constrained levels or
ground states.
\end{remark}

\subsection{Compactness and orbital stability}

The first lemma proves compactness of the ground states below the critical
mass.  The second establishes their orbital stability.  Then we combine these
results with Proposition~\ref{prop:cubic-level-classification} to prove
Theorem~\ref{thm:intro-cubic-ground-states}.

\begin{lemma}
\label{lem:cubic-ground-compactness}
Let \(0<\mu<\mu_*\).  Every minimizing sequence for
\(I_\alpha(\mu)\) has a subsequence converging strongly in \(X_\alpha\).
Consequently \(\mathcal M_\alpha(\mu)\) is nonempty and compact in
\(X_\alpha\).
\end{lemma}

\begin{proof}
Let \(M(u_n)=\mu\) and
\(E_\alpha(u_n)\to I_\alpha(\mu)\). \(u_n\) is bounded in \(X_\alpha\) by Lemma~\ref{lem:cubic-half-wave-coercivity}.  After
extraction,
\[
 u_n\rightharpoonup u\quad\text{in }X_\alpha,
\]
and write
\[
 u_n=\phi_n+q_nG_{2\lambda_\alpha},\qquad
 u=\phi+qG_{2\lambda_\alpha}.
\]
The norm identity \eqref{eq:critical-form-norm} gives
\[
 q_n\to q,\qquad
 \phi_n\rightharpoonup\phi\quad\text{in }H^{1/2}(\R).
\]
Lemma~\ref{lem:critical-form-coordinates} further gives
\[
 u_n\to u\quad\text{in }L^2(K)
\]
for every compact interval \(K\subset\R\).  A diagonal extraction therefore
also gives \(u_n\to u\) almost everywhere on \(\R\).

Suppose first that \(u=0\).  The preceding coordinate convergence gives
\(q_n\to0\), while
\(M(\phi_n)=\mu+o(1)<\mu_*\) for large \(n\).  Formula
\eqref{eq:critical-form} and the boundedness of \(\phi_n\) give
\[
 \mathfrak t_\alpha[u_n]
 =\norm{|D|^{1/2}\phi_n}_2^2+o(1),\qquad
 \norm{u_n}_4^4=\norm{\phi_n}_4^4+o(1).
\]
Therefore
\[
 E_\alpha(u_n)=E_\infty(\phi_n)+o(1)
 \geq I_\infty(M(\phi_n))+o(1)=o(1).
\]
This contradicts the negativity of \(I_\alpha(\mu)\) in
Proposition~\ref{prop:cubic-level-classification}.  Hence
\(m:=M(u)>0\).

Put \(r_n=u_n-u\).  Hilbert space splitting in the form norm, mass
splitting and the Brezis--Lieb lemma \cite{BrezisLieb} give
\[
 \mathfrak t_\alpha[u_n]
 =\mathfrak t_\alpha[u]+\mathfrak t_\alpha[r_n]+o(1).
\]
\[
 M(r_n)=\mu-m+o(1).
\]
\[
 E_\alpha(u_n)=E_\alpha(u)+E_\alpha(r_n)+o(1).
\]
Indeed, for the positive form inner product
\(\langle v,w\rangle_{X_\alpha}:=\mathfrak t_\alpha(v,w)
+2\lambda_\alpha\langle v,w\rangle\), weak convergence makes the cross
term \(\langle u,r_n\rangle_{X_\alpha}=o(1)\); subtracting the corresponding mass
splitting gives the first identity.  The last identity then follows from
Brezis--Lieb.

We now exclude dichotomy directly.  Suppose that \(0<m<\mu\), and write
\[
 r_n=\psi_n+s_nG_{2\lambda_\alpha},
 \qquad
 \psi_n:=\phi_n-\phi,
 \qquad
 s_n:=q_n-q\to0.
\]
Since \(s_n\to0\), \((\psi_n)\) is bounded in \(H^{1/2}\) and the fixed
function \(G_{2\lambda_\alpha}\) belongs to \(L^2\cap L^4\), the form, mass
and quartic errors satisfy
\[
 \mathfrak t_\alpha[r_n]-\norm{|D|^{1/2}\psi_n}_2^2=o(1),\qquad
 M(r_n)-M(\psi_n)=o(1),\qquad
 \norm{r_n}_4^4-\norm{\psi_n}_4^4=o(1).
\]
For the last relation, the standard quartic difference estimate gives
\[
 \left|\norm{r_n}_4^4-\norm{\psi_n}_4^4\right|
 \leq C|s_n|\norm{G_{2\lambda_\alpha}}_4
 \left(\norm{\psi_n}_4+|s_n|\norm{G_{2\lambda_\alpha}}_4\right)^3=o(1).
\]
Consequently,
\[
 E_\alpha(r_n)=E_\infty(\psi_n)+o(1),
 \qquad
 M(\psi_n)=\mu-m+o(1)<\mu_*
\]
for large \(n\).  Lemma~\ref{lem:free-cubic-geometry} therefore gives
\(\liminf_nE_\alpha(r_n)\geq0\).  The energy splitting implies
\[
 E_\alpha(u)\leq I_\alpha(\mu)
 <-\frac{\lambda_\alpha}{2}\mu.
\]
Set \(\vartheta:=\mu/m>1\).  Since
\[
 E_\alpha(\sqrt\vartheta u)-E_\alpha(u)
 =(\vartheta-1)
 \left(-\frac{\vartheta}{2}\mathfrak t_\alpha[u]
       +(\vartheta+1)E_\alpha(u)\right),
\]
the spectral bound \(\mathfrak t_\alpha[u]\geq-\lambda_\alpha m\) gives
\[
 -\frac{\vartheta}{2}\mathfrak t_\alpha[u]
 +(\vartheta+1)E_\alpha(u)
 <-\frac{\vartheta\lambda_\alpha\mu}{2}<0.
\]
Thus \(E_\alpha(\sqrt\vartheta u)<E_\alpha(u)\).  But
\(M(\sqrt\vartheta u)=\mu\), and hence
\[
 I_\alpha(\mu)
 \leq E_\alpha(\sqrt\vartheta u)
 <E_\alpha(u)
 \leq I_\alpha(\mu),
\]
a contradiction.  Thus \(m=\mu\).

Now \(M(r_n)\to0\), and the boundedness of \(r_n\) in \(X_\alpha\), together
with \(X_\alpha\hookrightarrow L^p\) for any fixed \(p>4\), gives
\(\norm{r_n}_4\leq\norm{r_n}_2^\theta\norm{r_n}_p^{1-\theta}\to0\), where
\(\theta=(p-4)/(2(p-2))\).
The spectral lower bound gives
\(\mathfrak t_\alpha[r_n]\geq-\lambda_\alpha\norm{r_n}_2^2=o(1)\), and hence
\(\liminf_nE_\alpha(r_n)\geq0\).  Taking the lower limit in the energy
splitting therefore shows
\[
 I_\alpha(\mu)\geq E_\alpha(u)\geq I_\alpha(\mu),
\]
so \(u\) is a minimizer.  The same splitting now gives
\(E_\alpha(r_n)\to0\); since \(\norm{r_n}_4\to0\), it follows that
\(\mathfrak t_\alpha[r_n]\to0\).  Finally,
\[
 \norm{r_n}_{X_\alpha}^2
 =\mathfrak t_\alpha[r_n]+2\lambda_\alpha\norm{r_n}_2^2\to0.
\]
This proves the asserted strong compactness.
\end{proof}

\begin{lemma}
\label{lem:cubic-ground-stability}
Let \(0<\mu<\mu_*\).  Then \(\mathcal M_\alpha(\mu)\) is orbitally stable
in \(X_\alpha\).
\end{lemma}

\begin{proof}
We use the Cazenave--Lions compactness argument \cite{CazenaveLions}.
First observe that solutions starting sufficiently close to
\(\mathcal M_\alpha(\mu)\) are global.  Indeed, since \(\mu<\mu_*\) and
\(X_\alpha\hookrightarrow L^2(\R)\), there is
\(\delta_0>0\) such that every \(u_0\in X_\alpha\) satisfying
\[
 \operatorname{dist}_{X_\alpha}
 \bigl(u_0,\mathcal M_\alpha(\mu)\bigr)<\delta_0
\]
also satisfies \(M(u_0)<\mu_*\).
Theorem~\ref{thm:intro-cubic-flow} then implies that the corresponding
solution is global.

Suppose now that \(\mathcal M_\alpha(\mu)\) is not orbitally stable.  Then
there exist \(\varepsilon_0>0\) and initial data \(u_{0,n}\in X_\alpha\) such
that
\[
 \operatorname{dist}_{X_\alpha}
 \bigl(u_{0,n},\mathcal M_\alpha(\mu)\bigr)\to0.
\]
After discarding finitely many terms, let \(u_n\) be the global solution with
\(u_n(0)=u_{0,n}\).  There are times \(t_n\in\R\) such that
\[
 \operatorname{dist}_{X_\alpha}
 \bigl(u_n(t_n),\mathcal M_\alpha(\mu)\bigr)\geq\varepsilon_0.
\]
By the definition of the distance, choose
\(\widetilde u_n\in\mathcal M_\alpha(\mu)\) such that
\[
 \norm{u_{0,n}-\widetilde u_n}_{X_\alpha}
 \leq \operatorname{dist}_{X_\alpha}
 \bigl(u_{0,n},\mathcal M_\alpha(\mu)\bigr)+\frac1n.
\]
Since \(M(\widetilde u_n)=\mu\) and
\(E_\alpha(\widetilde u_n)=I_\alpha(\mu)\), continuity of mass and energy on
\(X_\alpha\) yields
\[
 M(u_{0,n})\to\mu,
 \qquad
 E_\alpha(u_{0,n})\to I_\alpha(\mu).
\]
Set \(v_n:=u_n(t_n)\) and \(\mu_n:=M(v_n)\).  Conservation of mass and
energy gives
\[
 \mu_n=M(u_{0,n})\to\mu,
 \qquad
 E_\alpha(v_n)=E_\alpha(u_{0,n})\to I_\alpha(\mu).
\]
The masses \(M(v_n)\) converge to \(\mu\), which is strictly below \(\mu_*\).
Lemma~\ref{lem:cubic-half-wave-coercivity} therefore shows that \((v_n)\) is
bounded in \(X_\alpha\).  Define
\(a_n:=(\mu/\mu_n)^{1/2}\).  Then \(a_n\to1\) and
\[
 M(a_nv_n)=\mu.
\]
The boundedness of \((v_n)\) and the amplitude identity
\[
 E_\alpha(a_nv_n)-E_\alpha(v_n)
 =\frac{a_n^2-1}{2}\mathfrak t_\alpha[v_n]
  -\frac{a_n^4-1}{4}\norm{v_n}_4^4
\]
imply
\[
 E_\alpha(a_nv_n)\to I_\alpha(\mu).
\]
The same boundedness also gives
\[
 \norm{a_nv_n-v_n}_{X_\alpha}
 =|a_n-1|\norm{v_n}_{X_\alpha}\to0.
\]
Thus \((a_nv_n)\) is a minimizing sequence of exact mass \(\mu\).
Lemma~\ref{lem:cubic-ground-compactness} provides a subsequence and a
ground state \(\widetilde u\in\mathcal M_\alpha(\mu)\) such that
\[
 a_nv_n\to\widetilde u
 \quad\text{in }X_\alpha.
\]
For this subsequence, the triangle inequality gives
\[
 \operatorname{dist}_{X_\alpha}
 \bigl(v_n,\mathcal M_\alpha(\mu)\bigr)
 \leq\norm{v_n-a_nv_n}_{X_\alpha}
     +\norm{a_nv_n-\widetilde u}_{X_\alpha}
 \to0,
\]
which contradicts the defining property of \(t_n\).
\end{proof}

\begin{proof}[Proof of Theorem~\ref{thm:intro-cubic-ground-states}]
Let \(0<\mu<\mu_*\).  Lemma~\ref{lem:cubic-ground-compactness} shows that
every minimizing sequence is precompact in \(X_\alpha\), and hence
\(\mathcal M_\alpha(\mu)\) is nonempty.  Its orbital stability is
Lemma~\ref{lem:cubic-ground-stability}.  If \(\mu\geq\mu_*\),
Proposition~\ref{prop:cubic-level-classification} gives
\(I_\alpha(\mu)=-\infty\), so the constrained infimum cannot be attained and
\(\mathcal M_\alpha(\mu)=\varnothing\).
\end{proof}

\section{Near-critical rigidity toward free ground states}
\label{sec:critical-mass-rigidity}

This section proves the near-critical rigidity result.  We use a Green
decomposition adapted to the quadratic form and show that the rescaled
singular part vanishes.  The regular part then becomes an optimizing sequence
for the sharp free inequality, and compactness of the free optimizers gives
the result.

For \(v\in H^{1/2}(\R)\), define the free half-wave kinetic energy and the
quartic focusing interaction by
\[
 T(v):=\norm{|D|^{1/2}v}_2^2,
 \qquad
 N(v):=\norm v_4^4.
\]
Thus the free energy introduced in Section~\ref{sec:operator} is
\[
 E_\infty(v)=\frac12T(v)-\frac14N(v).
\]
The dilation \(\D_\rho\) preserves \(M\) and multiplies both \(T\) and
\(N\) by \(\rho\).

In addition to the sharp inequality \eqref{eq:intro-half-wave-GN}, we use the
Pohozaev normalization obtained as follows.  Testing the profile equation with
\(Q_*\) gives
\(T(Q_*)+M(Q_*)=N(Q_*)\).  Since \(Q_*\) is a critical point of
\(E_\infty\) under the fixed-mass constraint and \(\D_\rho\) preserves mass,
differentiation at \(\rho=1\) gives
\[
 \frac12T(Q_*)-\frac14N(Q_*)=0.
\]
Together with \(M(Q_*)=\mu_*\), these identities give
\begin{equation}\label{eq:rigidity-Q-pohozaev}
 T(Q_*)=M(Q_*)=\mu_*,
 \qquad
 N(Q_*)=2\mu_*.
\end{equation}

The following lemma gives the free compactness result needed below.

\begin{lemma}\label{lem:sharp-compactness}
Let \(v_n\) be bounded in \(H^{1/2}(\R)\) and assume
\begin{equation}\label{eq:optimizing-data}
 M(v_n)\to\mu_*,
 \qquad
 T(v_n)\to\mu_*,
 \qquad
 N(v_n)\to2\mu_*.
\end{equation}
Then
\begin{equation}\label{eq:compactness-general}
 \inf_{\theta,y\in\R}
 \norm{e^{-i\theta}v_n-Q_*(\,\cdot-y)}_{H^{1/2}}
 \to0.
\end{equation}
\end{lemma}

\begin{proof}
Apply the argument to an arbitrary subsequence.  Choose a smooth partition
\((\chi_j)_{j\in\mathbb Z}\) such that
\(\sum_j\chi_j^2=1\), each \(\chi_j\) is supported in \((j-1,j+1)\) and the
derivatives are uniformly bounded.  Then
\[
 \sum_j|\chi_j(x)-\chi_j(y)|^2
 \leq C\min\{1,|x-y|^2\},\qquad
 \sum_j\norm{\chi_jf}_{H^{1/2}}^2
 \leq C\norm f_{H^{1/2}}^2.
\]
Here the second estimate follows from the Gagliardo representation of the
\(H^{1/2}\)-seminorm.
Applying \eqref{eq:intro-half-wave-GN} to \(\chi_jf\), summing in \(j\) and
using finite overlap gives the concentration estimate
\[
 \norm f_4^4\leq C
 \left(\sup_{y\in\R}\norm f_{L^2(y-1,y+1)}^2\right)
 \norm f_{H^{1/2}}^2.
\]
Since \(N(v_n)\to2\mu_*>0\), translation and extraction give
\(v_n(\cdot+a_n)\rightharpoonup V\neq0\) in \(H^{1/2}\), with convergence
strongly in \(L^2_{\mathrm{loc}}\) and almost everywhere by local compactness.  For
\(r_n=v_n(\cdot+a_n)-V\), Hilbert space and Brezis--Lieb splitting
\cite{BrezisLieb} apply simultaneously to \(M,T,N\).  Put
\(M_1=M(V)\), \(T_1=T(V)\),
\(M_2=\mu_*-M_1\) and \(T_2=\mu_*-T_1\).  The sharp inequality gives
\[
 \mu_*^2\leq M_1T_1+M_2T_2\leq\mu_*^2.
\]
Indeed, with \(a=M_1/\mu_*\) and \(b=T_1/\mu_*\), the difference from the upper bound,
divided by \(\mu_*^2\), is \(a(1-b)+b(1-a)\).  Since \(V\neq0\), equality forces
\(a=b=1\).  Thus \(r_n\to0\) in \(H^{1/2}\), and \(V\) is an optimizer.
By \cite[Theorem~2.4]{FrankLenzmann}, every complex optimizer is a complex
multiple of a translate and dilation of \(Q_*\).  The identities
\(M(V)=T(V)=\mu_*\) fix the modulus and the dilation, so \(V\) is a phase
and translate of \(Q_*\).  Since the initial subsequence
was arbitrary, the full sequence approaches that orbit; choosing parameters
within \(1/n\) proves \eqref{eq:compactness-general}.
\end{proof}

The central mechanism of the rigidity argument is the vanishing of the
rescaled charge, which is established in the next lemma.

\begin{lemma}
\label{lem:near-critical-charge}
Let \(u_n\in X_\alpha\) satisfy the assumptions of
Theorem~\ref{thm:critical-rigidity}, and set
\(\tau_n=\mathfrak t_\alpha[u_n]\). For
\(u_n=\phi_n+q_nG_{\tau_n}\), we set
\(\rho_n=\mu_*/\tau_n\), \(\widetilde q_n=q_n/\sqrt{\tau_n}\),
\(L_n=\Phi_\alpha(-\tau_n)\) and \(w_n=\D_{\rho_n}\phi_n\).  Then
\begin{equation}\label{eq:near-charge}
 \abs{\widetilde q_n}=O\left(
 \frac1{L_n}+\sqrt{\frac{1-\mu_*/M(u_n)}{L_n}}
 \right),
\end{equation}
and \(w_n\) satisfy
\begin{equation}\label{eq:w-limits}
 M(w_n)\to\mu_*,
 \qquad
 T(w_n)\to\mu_*,
 \qquad
 N(w_n)\to2\mu_*.
\end{equation}
\end{lemma}

\begin{proof}
Since \(M(u_n)\to\mu_*\), the masses are bounded.  Hence
\eqref{eq:critical-form-norm} and
\(\norm{u_n}_{X_\alpha}\to\infty\) imply \(\tau_n\to\infty\).  In
particular, \(L_n>0\) for all sufficiently large \(n\).
Lemma~\ref{lem:green-coordinates} justifies the decomposition in
the statement with the same charge as in the coordinates associated with
\(G_{2\lambda_\alpha}\).  The form
identity
\eqref{eq:critical-form}, used with \(\lambda=\tau_n\), becomes
\begin{equation}\label{eq:near-form}
 1+M(u_n)=\frac{T(\phi_n)}{\tau_n}
 +M(\phi_n)+L_n\abs{\widetilde q_n}^2.
\end{equation}
Every term on the right is nonnegative.  Since \(M(u_n)\) is bounded,
\begin{equation}\label{eq:rough-r}
 \abs{\widetilde q_n}
 =O(L_n^{-1/2})=O((\log\tau_n)^{-1/2}).
\end{equation}

The Green kernels obey the exact scaling law
\begin{equation}\label{eq:green-scaling}
 G_\lambda(x)=G_1(\lambda x),
 \qquad
 \norm{G_\lambda}_s^s=\lambda^{-1}\norm{G_1}_s^s
 \quad(2\leq s<\infty).
\end{equation}
Put \(g_n=q_nG_{\tau_n}\).  It follows that
\begin{equation}\label{eq:g-norms}
 \norm{g_n}_2=\norm{G_1}_2\abs{\widetilde q_n},
 \qquad
 \norm{g_n}_4=\norm{G_1}_4\tau_n^{1/4}\abs{\widetilde q_n}.
\end{equation}
Since \(u_n=\phi_n+g_n\), the mass difference satisfies
\begin{equation}\label{eq:near-errors-mass}
 M(u_n)-M(\phi_n)=O(\abs{\widetilde q_n}).
\end{equation}
On the other hand, bounded energy and \eqref{eq:cubic-half-wave-energy} give
\begin{equation}\label{eq:u-N}
 N(u_n)=2\tau_n+O(1).
\end{equation}
The elementary quartic difference estimate, together with
\eqref{eq:g-norms} and \(\norm{u_n}_4=O(\tau_n^{1/4})\), yields
\begin{equation}\label{eq:N-error}
 \frac{N(\phi_n)}{\tau_n}
 =2+O(\abs{\widetilde q_n})+O(\tau_n^{-1}).
\end{equation}

We now feed these estimates into the sharp free inequality.  Equations
\eqref{eq:intro-half-wave-GN}, \eqref{eq:near-errors-mass} and
\eqref{eq:N-error} imply
\begin{equation}\label{eq:near-T}
 \frac{T(\phi_n)}{\tau_n}
 \geq\frac{\mu_*}{M(u_n)}-C\abs{\widetilde q_n}-C\tau_n^{-1}.
\end{equation}
Using the lower estimates \eqref{eq:near-errors-mass} and \eqref{eq:near-T}
in \eqref{eq:near-form} gives
\begin{equation}\label{eq:near-quadratic}
 L_n\abs{\widetilde q_n}^2
 \leq1-\frac{\mu_*}{M(u_n)}+C\abs{\widetilde q_n}+C\tau_n^{-1}.
\end{equation}
Solving this quadratic inequality first gives
\[
 \abs{\widetilde q_n}
 \lesssim \frac1{L_n}
 +\sqrt{\frac{1-\mu_*/M(u_n)}{L_n}}
 +(\tau_nL_n)^{-1/2}.
\]
Since \(L_n=O(\log\tau_n)=o(\tau_n)\), the last term is
\(O(L_n^{-1})\).  This proves \eqref{eq:near-charge}.
Since \(M(u_n)\to\mu_*\), equations \eqref{eq:near-charge} and
\eqref{eq:near-quadratic} give \(\widetilde q_n\to0\) and
\(L_n\abs{\widetilde q_n}^2\to0\).

Equations \eqref{eq:near-errors-mass}, \eqref{eq:N-error} and
\eqref{eq:near-form} now give
\begin{equation}\label{eq:normalized-limits}
 M(\phi_n)\to\mu_*,
 \qquad
 \frac{T(\phi_n)}{\tau_n}\to1,
 \qquad
 \frac{N(\phi_n)}{\tau_n}\to2.
\end{equation}
Applying the scaling laws to \(w_n=\D_{\rho_n}\phi_n\) proves
\eqref{eq:w-limits}.
\end{proof}

\begin{proof}[Proof of Theorem~\ref{thm:critical-rigidity}]
Apply Lemma~\ref{lem:near-critical-charge}.  The limits
\eqref{eq:w-limits} and Lemma~\ref{lem:sharp-compactness} give
\begin{equation}\label{eq:regular-convergence}
 \inf_{\theta,y\in\R}
 \norm{e^{-i\theta}\D_{\rho_n}\phi_n-Q_*(\,\cdot-y)}_{H^{1/2}}
 \to0.
\end{equation}

Lemma~\ref{lem:green-coordinates} gives
\begin{equation}\label{eq:scaled-green}
 \D_{\rho_n}(q_nG_{\tau_n})
 =\sqrt{\mu_*}\,\widetilde q_nG_{\mu_*}.
\end{equation}
Here the dilation is used only to normalize a static sequence and to send the
matched Green kernel to \(G_{\mu_*}\).  It is not being
treated as a symmetry of the operator with fixed coupling \(\alpha\).
Since \(\widetilde q_n\to0\), the right-hand side tends to zero in the form
norm.  Indeed,
Lemma~\ref{lem:green-coordinates} gives
\(G_{\mu_*}-G_{2\lambda_\alpha}\in H^{1/2}(\R)\).  Thus, in the
coordinates associated with \(G_{2\lambda_\alpha}\), the additional regular
part tends to zero in \(H^{1/2}(\R)\), while its charge tends to zero.
The explicit norm formula \eqref{eq:critical-form-norm}, together with
\eqref{eq:regular-convergence}, proves \eqref{eq:full-X}.
\end{proof}

\begin{remark}
The condition \(M(u_n)\to\mu_*\) in Theorem~\ref{thm:critical-rigidity} is
essential.  Fix a mass \(m>\mu_*\).  One can choose an even
\(V\in H^{1/2}(\R)\) with \(M(V)=m\) and \(E_\infty(V)=0\) which is not a
free optimizer.  Choose an even integer \(N\) such that \(m/N<\mu_*\), and
place \(N\) nonnegative smooth compactly supported packets of mass \(m/N\) in
symmetric pairs with mutually large separations.  Their sum \(V_+\) is even
and has mass \(m\).  The sharp free inequality makes the energy of each packet
positive.  The quartic terms are additive for disjoint supports, while the
kinetic cross terms vanish as the separations tend to infinity.  Indeed, the
Fourier representation and the Riemann--Lebesgue lemma give the latter
convergence.  Thus the
packets can be chosen so that \(E_\infty(V_+)>0\).  On the other hand, the even
function
\[
 V_-:=\sqrt{m/\mu_*}\,Q_*
\]
has mass \(m\) and
\[
 E_\infty(V_-)=\frac m2\left(1-\frac m{\mu_*}\right)<0.
\]
Both functions are nonnegative.  Hence the denominator in the continuous
fixed-mass path
\[
 \Gamma(\theta):=\sqrt m\,
 \frac{(1-\theta)V_++\theta V_-}
 {\norm{(1-\theta)V_++\theta V_-}_2},
 \qquad 0\leq\theta\leq1,
\]
does not vanish.  The continuity of \(E_\infty\) on \(H^{1/2}(\R)\) and the
intermediate value theorem give \(V=\Gamma(\theta_0)\) with
\(E_\infty(V)=0\).  It cannot be a free optimizer: equality in the sharp free
inequality would instead imply
\[
 E_\infty(V)=\frac12\left(1-\frac m{\mu_*}\right)T(V)<0.
\]

Let \(s_n\to\infty\) and set \(u_n=\D_{s_n}V\).  These functions have zero
charge, and hence
\[
 M(u_n)=m,\qquad E_\alpha(u_n)=0,
 \qquad \mathfrak t_\alpha[u_n]=s_nT(V)\to\infty.
\]
In particular, \(\norm{u_n}_{X_\alpha}\to\infty\).  Their form normalization
is the fixed profile
\[
 \D_{\mu_*/\mathfrak t_\alpha[u_n]}u_n
 =W:=\D_{\mu_*/T(V)}V.
\]
Since \(E_\infty(V)=0\), its sharp inequality deficit is
\[
 \frac{2M(W)}{\mu_*}T(W)-N(W)=2(m-\mu_*)>0.
\]
Thus \(W\) is not a free optimizer, and the rigidity conclusion does not
extend to any fixed supercritical mass.
\end{remark}

\end{document}